\documentclass[final]{mcom-l}

\usepackage{amssymb}
\usepackage{mythesis}
\newcommand{\ba}{\begin{array}}
\newcommand{\ea}{\end{array}}
\newcommand{\bd}{\begin{displaymath}}
\newcommand{\ed}{\end{displaymath}}
\newcommand{\bi}{\begin{itemize}}
\newcommand{\ei}{\end{itemize}}
\newcommand{\bn}{\begin{enumerate}}
\newcommand{\en}{\end{enumerate}}
\newcommand{\pa}{\partial}

\newcommand{\f}{\frac}

\UseRawInputEncoding

\usepackage{graphicx}
\usepackage{float}
\newtheorem{theorem}{Theorem}[section]
\newtheorem{lemma}[theorem]{Lemma}
\newtheorem{corollary}[theorem]{Corollary}

\theoremstyle{definition}

\newtheorem{prob}{Problem}

\theoremstyle{remark}
\newtheorem{remark}[theorem]{Remark}

\numberwithin{equation}{section}

\begin{document}

\title[Analysis of Iterative Solution of Helmholtz via the Wave Equation]{Extensions and Analysis of an Iterative Solution of the Helmholtz Equation via the Wave Equation}

\author{Fortino Garcia}
\address{Department of Applied Mathematics, University of Colorado, Boulder, CO 80309-0526}
\email{Fortino.Garcia@colorado.edu}
\thanks{The first author is supported in part by STINT initiation grant IB2019--8154 and NSF Grant DMS-1913076. Any conclusions or recommendations expressed in this paper are those of the author and do not necessarily reflect the views of the NSF.}

\author{Daniel Appel\"o}
\address{Department of Computational Mathematics, Science, and Engineering; Department of Mathematics, Michigan State University, East Lansing, MI 48824}
\email{appeloda@msu.edu}
\thanks{The second author is supported in part by STINT initiation grant IB2019--8154 and NSF Grant DMS-1913076. Any conclusions or recommendations expressed in this paper are those of the author and do not necessarily reflect the views of the NSF.}

\author{Olof Runborg}
\address{Department of Mathematics, KTH, 100 44, Stockholm}
\email{olofr@kth.se}
\thanks{The third author is supported in part by STINT initiation grant IB2019--8154.}




\begin{abstract}
In this paper we extend analysis of the WaveHoltz iteration --
a time-domain iterative method for the solution of the Helmholtz equation.
We expand the previous analysis of energy conserving problems
and prove convergence of the WaveHoltz iteration for problems with impedance 
boundary conditions in a single spatial dimension. We then consider 
interior Dirichlet/Neumann problems with
damping in any spatial dimension, and show that 
for a sufficient level of damping the WaveHoltz iteration converges in a number 
of iteration independent of the frequency. Finally, we present a discrete analysis 
of the WaveHoltz iteration for a family of higher order time-stepping schemes.
We show that the fixed-point of the discrete WaveHoltz iteration 
converges to the discrete Helmholtz solution with the order of the 
time-stepper chosen. We present numerical examples and demonstrate that it is possible to 
\textit{completely remove} time discretization error from the 
WaveHoltz solution through careful analysis of the discrete iteration
together with updated quadrature formulas. 
\end{abstract}

\maketitle	


\section{Introduction}
This is the second of a series of papers on time-domain 
methods for the numerical solution of the Helmholtz equation 
\begin{equation}
  \nabla\cdot(c^2(x)\nabla u) + \omega^2u = f(x),\qquad
  x\in \Omega, \label{eq:helm}
\end{equation}
for a domain $\Omega$, frequency $\omega$, and sound speed $c(x)$. 
The Helmholtz equation (both acoustic and elastic) is useful for seismic, 
acoustic, and optics applications. The numerical solution of 
the Helmholtz equation is especially difficult due to 
the resolution requirements and the indefinite nature of the 
Helmholtz operator for large frequencies. 

In the previous paper \cite{WaveHoltz}, we introduced a time-domain 
approach for solving the Helmholtz equation \eqref{eq:helm}.
Given the Helmholtz solution, $u(x)$, 
the time-harmonic wave field $\text{Re}\{u(x)e^{-i\omega t}\}$
satisfies the wave equation
\begin{eqnarray}
&&  w_{tt} =   \nabla\cdot(c^2(x)\nabla w)
  - f(x) \cos(\omega t), \quad x\in\Omega, \ \ 0 \le t \le T, \nonumber \\
&&  w(0,x) = v_0(x), \quad w_t(0,x) = v_1(x), \nonumber
\end{eqnarray}
where $v_0 = \text{Re}\{u(x)\}$ and $v_1 = \omega\text{Im}\{u(x)\}$.
In \cite{WaveHoltz}, we introduced an integral operator that time-filtered
the wave solution resulting from initial data $v^n_0$, $v^n_1$. The time-filtering
generates new iterates $v^{n+1}_0$, $v^{n+1}_1$ leading to a fixed-point iteration 
we named the WaveHoltz iteration. The convergence of the fixed-point iteration 
for interior problems with Dirichlet/Neumann boundary conditions (i.e. energy conserving problems)
was proven in the continuous and discrete settings
under a non-resonance condition. For such problems, the WaveHoltz 
iteration can be reformulated as a symmetric and 
positive-definite system which can be accelerated with Krylov subspace
methods such as the conjugate gradient method and GMRES. Numerical experiments 
using the WaveHoltz iteration indicated promising scaling with frequency
for problems also with outflow boundary conditions common in seismic applications,
though no theoretical proof was given for the convergence of the method 
in that case.

In this paper, we extend the continuous and discrete analysis 
presented in the prequel.
In~\cite{WaveHoltz}, the continuous analysis was performed using 
a simplified iteration in which the initial velocity, $v^n_{1}$,
is set to zero each iteration as solutions to the Helmholtz equation are
real-valued. In this paper we begin by proving convergence of the 
WaveHoltz iteration without the assumption $v^n_{1} \equiv 0$.
This result, together with appropriate extensions of the problem 
data, leads to a proof of convergence of the WaveHoltz iteration
for problems with impedance boundary conditions in a single spatial dimension.
To conclude the continuous analysis, we additionally consider 
the damped Helmholtz equation and prove that the iteration is 
convergent. Numerical results verify that for a sufficiently large damping, the 
number of iterations for the WaveHoltz iteration to reach convergence for damped 
Helmholtz equations is independent of frequency. We thus can guarantee convergence of the method
to the Helmholtz solution via impedance conditions and/or damping,
without any additional conditions, 
which we note is possible due to the absence of resonant frequencies.

For the discrete analysis, we investigate the effect of choice of time-stepper 
used for the WaveHoltz iteration. In~\cite{WaveHoltz}, we noted that 
in the discrete case the WaveHoltz iteration converged to the solution
of a discrete Helmholtz problem with modified frequency. We provided 
the modification for a centered second order time-stepping scheme which would 
recover the original discrete Helmholtz solution. Here
we consider higher order modified equation (ME) time-stepping schemes, 
\cite{shubin1987modified,anne2000construction},
and show that the fixed-point of the discrete WaveHoltz iteration 
converges to the discrete Helmholtz solution with the order of the 
time-stepper chosen. We additionally show that, 
as in the case for EM-WaveHoltz~\cite{peng2021emwaveholtz}, 
it is possible to \textit{completely remove} time discretization error from the 
WaveHoltz solution through careful analysis of the discrete iteration
and updated quadrature formulas.

The efficient solution of the Helmholtz equation \eqref{eq:helm}
via iterative methods is notoriously difficult, especially for
high-frequency problems of practical interest, and has been 
the subject of much research. We refer to our previous paper \cite{WaveHoltz} 
for a more in-depth overview of the literature on techniques for 
solving the Helmholtz equation, as well as the 
review articles \cite{ernst2012difficult,Gander_Zhang_SIAM_REV,erlangga2008advances}.
We focus on the literature that is closely related to the methods and approach used here.

The theoretical justification for working in the time-domain comes from 
the \textit{limiting amplitude principle}, see 
\cite{morawetz1962limiting,Ladyzhenskaya:57,Vainberg:75}.
The principle states that every solution to the wave equation 
with a time-harmonic forcing in the exterior of a domain with reflecting boundary conditions
tends to the Helmholtz solution. 
Rather than evolving a wave equation forward in time to reach a steady state by 
appealing to the limiting amplitude principle, it is possible to 
cast the problem as a constrained convex least-squares minimization problem.
This approach, originally proposed by Bristeau et al. \cite{bristeau1998controllability},
is the so-called Controllability Method (CM). The CM seeks to accelerate the convergence  
to the steady-state limit by minimizing the deviation from
time-periodicity of the time-domain solution in second-order form. 

In the original CM, along with later work by Heikkola et al. 
\cite{HEIKKOLA2007344,HEIKKOLA20071553}, only scatterers with Dirichlet boundary conditions 
were considered as the original cost functional 
of \cite{bristeau1998controllability} did not generally yield unique 
minimizers for other types of boundary conditions. An alternative 
functional, $J_\infty$, proposed by Bardos and Rauch in \cite{bardos1994variational}, 
however, did yield uniqueness of the minimizer at the cost of requiring the storage 
of the entire history of the computed solution to the wave equation 
which could be prohibitive for large problems.

For the wave equation in second-order form, the initial condition lies in 
$H^1 \times L^2$, requiring the solution of a coercive elliptic problem 
to find a Riesz representative for gradient calculations. Glowinski and 
Rossi \cite{GlowRoss06} presented an update to the CM by considering the wave 
equation in first-order form, allowing the initial conditions to lie in a reflexive space
and thus removing the need for an elliptic solve in each iteration.
The discretization chosen in this case, however, had the drawback of 
requiring inversion of a mass-matrix at each time-step.

In more recent work by Grote and Tang, \cite{grote2019controllability}, the use
of an alternative functional (or post-processing via a compatibility condition)
restored uniqueness of the minimizer of CM. In a follow-up paper, \cite{GROTE2020112846}, Grote et al. 
proposed the use of a hybrid discontinuous Galerkin discretization, \cite{Stanglmeier2016}, of the 
first-order form wave equation which allowed the scheme to be fully explicit and therefore 
fully parallel. Moreover, they extend CM to general boundary conditions for the first-order
formulation and additionally proposed a filtering procedure which allows the 
original energy functional to be used regardless of the boundary condition.

The above work has inspired other time-domain methods outside of CM and WaveHoltz.
Work by Stolk \cite{stolk2020timedomain} utilizes
time-domain approaches as a preconditioner for a GMRES accelerated preconditioner for
direct Helmholtz discretizations yielding a  hybrid time-frequency domain method. 
Arnold et al. \cite{arnold2021adaptive} propose a time-domain method for scattering
problems which leverages the compact support of incident field plane wavelets together with 
a front-tracking adaptive meshing algorithm to reduce the cost of computing a 
Fourier transform of the wave solution to obtain Helmholtz solutions.

Another important class of methods for solving the Helmholtz equation are
the so-called shifted Laplacian preconditioners. The use of the Laplacian 
as a preconditioner for Helmholtz problems emerged with the initial work of
Bayliss et al. \cite{bayliss1983iterative}. In \cite{bayliss1983iterative},
the normal equations of the discrete Helmholtz equation were iteratively solved
using conjugate gradient, with a Symmetric Successive Over-Relaxation (SSOR) 
sweep of the discrete Laplacian as a preconditioner. Giles and Laird then
extended the previous preconditioner to instead solve the Helmholtz system 
with a flipped sign in front of the Helmholtz term using multigrid \cite{laird2002preconditioned}. 
Erlangga, Vuik and Osterlee \cite{erlangga2004class,Erlangga2006} further generalized the
previous work to use a complex-valued shift of the Laplacian leading to the
shifted Laplacian preconditioner. For a review of the class of shifted Laplacian
preconditioners we refer the reader to the review article by Erlangga \cite{erlangga2008advances}.

The rest of this paper is organized as follows. In Section 2 we present analysis for the 
general WaveHoltz iteration and prove convergence in the case of impedance boundary
conditions in a single spatial dimension. In Section 3 we present a brief analysis for the case in which 
damping is present. Section 4 outlines a discrete analysis of higher order
modified equation (ME) schemes, and we additionally present a method to \textit{completely}
remove time discretization error from the discrete WaveHoltz solution. Finally, in Section 5 
we describe our numerical methods, Section 6 present our numerical examples, and summarize 
the paper in Section 7.

\section{The General Iteration}\label{sec::GeneralIteration}
We consider the Helmholtz equation in a bounded open smooth
domain $\Omega$,
\begin{equation}
  \nabla\cdot(c^2(x)\nabla u) + \omega^2u = f(x),\qquad
  x\in \Omega, \label{eq::helm}
\end{equation}
with boundary conditions of the type
\begin{equation}
i\alpha\omega u + \beta(c(x)\vec{n}\cdot \nabla u)=0,
\quad \alpha^2+\beta^2=1,\quad
\qquad x\in \partial\Omega.
 \label{eq::helmbc}
\end{equation}
We assume $f\in L^2(\Omega)$ and that $c\in L^\infty(\Omega)$ with
the bounds $0< c_{\min}\leq c(x)\leq c_{\max}<\infty$ a.e. in $\Omega$.
Away from resonances, this ensures that there is 
a unique weak solution $u\in H^1(\Omega)$
to \eqref{eq::helm}.
Due to the boundary conditions $u$ is in general complex-valued.

We first note that the function $w(t,x) := \text{Re}\{u(x) \exp(-i\omega t)\}$
is a $T=2\pi/\omega$-periodic (in time) 
solution to the real-valued forced scalar wave equation
\begin{eqnarray}
&&  w_{tt} =   \nabla\cdot(c^2(x)\nabla w)
  - \text{Re}\{f(x) e^{-i\omega t}\}, \quad x\in\Omega, \ \ 0 \le t \le T, \nonumber \\
&&  w(0,x) = v_0(x), \quad w_t(0,x) = v_1(x), \nonumber \\
&&  \alpha w_t +\beta(c(x)\vec{n}\cdot \nabla w)=0, \quad x\in\partial\Omega
, \label{eq:wave}
\end{eqnarray}
where $v_0=\text{Re}\{u\}$ and $v_1=\omega \text{Im}\{u\}$.
Based on this observation, our approach is to find this $w$
instead of $u$. We could thus look for
initial data  $v_0$ and $v_1$ such that $w$ is a $T$-periodic solution
to (\ref{eq:wave}). 
However, there may be several such $w$, see \cite{grote2019controllability}, 
and we therefore impose
the alternative constraint that a certain time-average of $w$
should equal the initial data. 
More precisely, we introduce the
following operator acting on the initial data $v_0\in H^1(\Omega)$, $v_1\in L^2(\Omega)$,
\be\label{eq::filterstep}
{\Pi} \left[
\begin{array}{c}
v_0\\
v_1 
\end{array}
\right] = \frac{2}{T}\int_0^{T} \left(\cos(\omega t) - \frac{1}{4}\right) \left[
\begin{array}{c}
w(t,x)\\
w_t(t,x) 
\end{array}
\right] dt,\quad
T=\frac{2\pi}{\omega},
\ee
where $w(t,x)$ and its time derivative $w_t(t,x)$ satisfies the wave equation (\ref{eq:wave}) with initial data $v_0$ and $v_1$.
The result of $\Pi [v_0,\ v_1]^T$ can thus be seen as
a filtering in time 
of $w(\cdot,x)$
around the $\omega$-frequency.
By construction, the solution $u$ of Helmholtz now satisfies the
system of equations
\begin{equation}
\left[
\begin{array}{c}
\text{Re}\{u\}\\
\omega \text{Im}\{u\} 
\end{array}
\right] = {\Pi} \left[
\begin{array}{c}
\text{Re}\{u\}\\
\omega \text{Im}\{u\} 
\end{array}
\right].
\end{equation}
The WaveHoltz iteration then amounts to solving
this system of equations with the fixed point iteration
\begin{equation}\label{eqn::Iteration}
\boxed{
\left[
\begin{array}{c}
v\\
v' 
\end{array}
\right]^{(n+1)} = {\Pi} \left[
\begin{array}{c}
v\\
v' 
\end{array}
\right]^{(n)}, \qquad \left[
\begin{array}{c}
v\\
v' 
\end{array}
\right]^{(0)} \equiv 0.
}
\end{equation}
Provided this iteration converges and the solution to
is unique, we 
obtain the Helmholtz solution as
\mbox{$u=\lim_{n\to\infty} v^n$}.

As seen in \cite{WaveHoltz}, the WaveHoltz operator $\Pi$ is affine
and can be written as $\Pi\v = {\Scal}\v+b$, where $\Scal$
is a linear operator, $\v=(v, v')^T$ and $b$ a fixed function.
Since the sought solution satisfies $\Pi\v=\v$, we can then reformulate the iteration
as a linear system 
\begin{align*}
    \mathcal{A}\v = b, \qquad \mathcal{A}=I-\Scal,
  \end{align*}
which allows the convergence of the WaveHoltz iteration 
to be accelerated by a Krylov method.
We note that the right hand side can be computed by applying
the WaveHoltz operator to the zero function, $b=\Pi 0$.
The action of $\mathcal{A}$ can also be computed via one
application of $\Pi$, as
 $\mathcal{A}\v = \v-\Pi\v+b$. Hence, after precomputing $b$
 the action can be computed by applying $\Pi$ to $\v$,
 i.e. by evolving the
wave equation for one period in time with initial data $\v$ and filter
the solution.
There is no
need to explicitly form $\mathcal{A}$. 

\begin{remark}
The operator $\mathcal{A}$ for the general iteration is not 
symmetric unlike the simplified iteration for energy conserving
problems where $v_1 = 0$. For interior, energy conserving problems
we recommend the use of the simplified iteration so that 
the conjugate gradient method may be used to accelerate convergence. 
For other boundary conditions, the general 
WaveHoltz iteration is required and a more versatile Krylov method,
such as GMRES, should be used.
\end{remark}

\subsection{Convergence for the Energy Conserving Case for the General WaveHoltz Iteration} \label{sec:iteration}

Here we consider boundary conditions of either
Dirichlet ($\beta=0$) or Neumann ($\alpha=0$) type
in \eqref{eq:wave}. 
This is typically
the most difficult case for iterative Helmholtz solvers
when $\Omega$ is bounded.
The wave energy is preserved in time and certain $\omega$-frequencies
in Helmholtz are resonant, meaning they equal an eigenvalue
of the operator $-\nabla\cdot(c^2(x)\nabla)$. Moreover,
the limiting amplitude principle does not hold, and one
can thus not obtain the Helmholtz solution by solving the
wave equation over a long time interval. We note that convergence 
of the WaveHoltz iteration in the energy conserving case was proved in~\cite{WaveHoltz}
using a simplified iteration for which $v'^{(n)} \equiv 0$ in \eqref{eqn::Iteration}.
In this section we prove convergence of the general iteration \eqref{eqn::Iteration}
without the assumption that $v'^{(n)} \equiv 0$. With this result in hand,
it will then be possible to establish convergence for the non-energy conserving case 
in Section~\ref{sec::non-energy}.

By the choice of boundary conditions the operator 
$-\nabla\cdot(c^2(x)\nabla)$ has a point spectrum with non-negative
eigenvalues. Denote those eigenmodes $(\lambda_j^2,\phi_j(x))$. 
We assume that the angular frequency $\omega$ is not a resonance, i.e. $\omega^2\neq \lambda_j^2$ for all $j$. 
The Helmholtz equation (\ref{eq:helm}) is then wellposed. 

We recall that for any $q\in L^2(\Omega)$ we can
expand
$$
  q(x) = \sum_{j=0}^\infty q_j\phi_j(x),
$$
for some coefficients $q_j$ and
$$
   ||q||_{L^2(\Omega)}^2 = \sum_{j=0}^\infty |q_j|^2,\qquad
   c_{\min}^2||\nabla q||_{L^2(\Omega)}^2
\leq   \sum_{j=0}^\infty \lambda_j^2|q_j|^2\leq 
      c_{\max}^2||\nabla q||_{L^2(\Omega)}^2.
$$
We start by expanding the Helmholtz solution $u = u^R + i u^I$, the initial data $v_0, \, v_1$ 
to the wave equation (\ref{eq:wave}), and the forcing $f = f^R + i f^I$ in this way,
\begin{gather*}
  u^R(x) = \sum_{j=0}^\infty u^R_j \phi_j(x), \ \ \ \  v_0(x) = \sum_{j=0}^\infty v_{0,j} \phi_j(x), \\
  v_1(x) = \sum_{j=0}^\infty v_{1,j} \phi_j(x), \ \ \ \    f^R(x) = \sum_{j=0}^\infty f^R_j \phi_j(x),
\end{gather*}
with analogous expansions for the imaginary parts of $u$ and $f$, $u^I$ and 
$f^I$, respectively.
Then,
\[
-\lambda_j^2u^R_j+\omega^2u^R_j = f^R_j \quad\Rightarrow\quad
u^R_j = \frac{f^R_j}{\omega^2-\lambda_j^2},
\]
and similarly for the imaginary parts $u^I_j$ and $f^I_j$.
For the wave equation solution $w(t,x)$ with initial data $w=v_0$ and $w_t=v_1$ we have
\begin{gather*}
  w(t,x)  = \sum_{j=0}^\infty w_j(t) \phi_j(x),
\end{gather*}
where
\begin{align*}
  w_j(t) 
   =u^R_j \left[\cos(\omega t) -\cos(\lambda_j t) \right]
  &+ u^I_j \left[\sin(\omega t) -\frac{\omega}{\lambda_j}\sin(\lambda_j t)\right]
  \\
  &+ v_{0,j}\cos(\lambda_j t)
  + \frac{v_{1,j}}{\lambda_j}\sin(\lambda_j t), \label{eq:no_filter_yet}
\end{align*}
with
  \begin{align*}
  w_0^\text{Neu}(t) 
   =u^R_0 \left[\cos(\omega t) -1 \right]
  + u^I_0 \left[\sin(\omega t) -\omega t\right]
  + v_{0,0}
  + v_{1,0}t,
  \end{align*}
if $\lambda_0 = 0$, as is the case for Neumann boundary conditions (a special 
case which we denote via the superscript `$\text{Neu}$' in the following analysis).
The filtering step then gives
\begin{gather*}
  \Pi \begin{bmatrix}
    v_0 \\ v_1
  \end{bmatrix}
  =
  \sum_{j=0}^\infty
  \begin{bmatrix}
    \bar{v}_{j}  \\
    \bar{v}'_{j}
  \end{bmatrix}
  \phi_j(x),
\end{gather*}
where
\begin{align*}
    \bar{v}_{j} & = u^R_j\left(1-\beta(\lambda_j) \right) - u^I_j \frac{\omega}{\lambda_j} \gamma(\lambda_j)+ v_{0,j}\beta(\lambda_j) +  \frac{v_{1,j}}{\lambda_j} \gamma(\lambda_j), \\
    \bar{v}'_{j} & = u^R_j \lambda_j\gamma(\lambda_j) + \omega u^I_j\left(1 - \beta(\lambda_j)\right) - v_{0,j}\lambda_j\gamma(\lambda_j) +  v_{1,j} \beta(\lambda_j),
\end{align*}
and
\[
\beta(\lambda)
:= \frac{2}{T}\int_0^T \left(\cos(\omega t) - \frac{1}{4}\right)\cos(\lambda t) dt, \quad 
\gamma(\lambda)
:= \frac{2}{T}\int_0^T \left(\cos(\omega t) - \frac{1}{4}\right)\sin(\lambda t) dt.
 \]
when $\lambda_j \ne 0$. For the Neumann case, we have
\begin{align*}
    \bar{v}_{0} = \frac{3}{2}u^R_0 + \frac{\pi}{2}u^I_0 - \frac12 v_{0,0} -\frac{\pi}{2\omega}v_{1,0}, \qquad
    \bar{v}'_{0} = \frac{3}{2} u^I_j - \frac12 v_{1,j}.
\end{align*}
By definition we have 
  \begin{align}\label{eqn::gammaBound}
    \left|\frac{\gamma(\lambda_j)}{\lambda_j}\right| \le 
    \frac{2}{T}\int_0^T \left|\left(\cos(\omega t) - \frac{1}{4}\right)\right|\left|t\frac{\sin(\lambda t)}{\lambda t}\right| dt
    \le 
   \frac{2}{T}\int_0^T \frac{5}{4} t dt
   =
   \frac{5 \pi}{2\omega},
  \end{align}
since $|\sin(x)/x|\le1$, which ensures the boundedness of the coefficients $\bar{v}_{j}, \bar{v}'_{j}$ 
for small eigenvalues $\lambda_j$.

Letting $v_{0,j}, v_{1,j}$ denote the coefficients of $v_{0}, v_{1}$ in the eigenbasis of the Laplacian, we can write the iteration as
\begin{align}\label{eqn::IterationMode}
  \begin{bmatrix}
    v_{0,j}^{n+1}\\v_{1,j}^{n+1}
  \end{bmatrix}
  =
  \left(\Pi
  \begin{bmatrix}
    v_{0}^{n}\\v_{1}^{n}
  \end{bmatrix}\right)_j
  =
  \left(I - B_j\right)
  \begin{bmatrix}
    u^R_j\\\omega u^I_j
  \end{bmatrix}
  +
  B_j
  \begin{bmatrix}
    v_{0,j}^{n}\\v_{1,j}^{n}
  \end{bmatrix},
\end{align}
where if we define $\beta_j = \beta(\lambda_j)$ and $\gamma_j = \gamma(\lambda_j)$ then
\begin{align*}
   B_j =  
   \begin{pmatrix}
    \beta_j & \gamma_j/\lambda_j \\
    -\lambda_j\gamma_j & \beta_j \\
   \end{pmatrix},
   \quad 
   B_0^\text{Neu} = 
   \begin{pmatrix}
    -1/2 & -\pi/2\omega \\
    0 & -1/2 \\
   \end{pmatrix},
\end{align*}
Moreover, the eigenvectors and eigenvalues of $B_j$
are
  \begin{align*}
    \xi_j^{\pm} = \begin{pmatrix}
      \pm i/\lambda \\ 1
    \end{pmatrix},
    \quad 
    \xi_0^\text{Neu} = \begin{pmatrix}
      1 \\ 0
    \end{pmatrix},
    \quad 
    \mu_j = \beta_j \pm i \gamma_j.
  \end{align*}
Introducing the 
linear operator ${\mathcal S}:L^2(\Omega)\times L^2(\Omega) \to L^2(\Omega)\times L^2(\Omega)$,
\begin{align}\label{eqn::SpectralOp}
  \mathcal{S} \sum_{j=0}^\infty \begin{bmatrix}  u^R_j \\ u^I_j  \end{bmatrix}\phi_j(x)
  =  \sum_{j=0}^\infty B_j \begin{bmatrix} u^R_j \\  u^I_j  \end{bmatrix}\phi_j(x),
\end{align}
we may write the iteration as
\begin{equation}\label{eqn:iteration}
  \boxed{
  \left[
  \begin{array}{c}
  v\\
  v' 
  \end{array}
  \right]^{(n+1)} = {\Pi} \left[
  \begin{array}{c}
  v\\
  v' 
  \end{array}
  \right]^{(n)}
  =
  \left[
  \begin{array}{c}
  u^R\\
  \omega u^I 
  \end{array}
  \right]
  +
  \mathcal{S}\left(
  \left[
  \begin{array}{c}
  v\\
  v' 
  \end{array}
  \right]^{(n)}
  -
  \left[
  \begin{array}{c}
  u^R\\
  \omega u^I 
  \end{array}
  \right]
  \right).
}
\end{equation}
We note that, in contrast to the simplified iteration 
analyzed in \cite{WaveHoltz},
the operator $\mathcal{S}$ is not symmetric for the 
general iteration. Despite this, we may identify the eigenmodes 
of $\mathcal{S}$ from the eigenvectors
of $B_j$ via $ \xi_j^{\pm} \phi_j$ with eigenvalues 
$\mu_j =\beta_j \pm i \gamma_j$
and $\xi_0^\text{Neu} = \xi_0^\text{Neu} \phi_0$ with eigenvalue
$\mu_0^\text{Neu} = -1/2$.

From \eqref{eqn::IterationMode}, we see 
that the iteration for each mode takes the form 
\begin{align*}
  \begin{bmatrix}
    v_{0,j}^{n+1}\\v_{1,j}^{n+1}
  \end{bmatrix}
  =
  \left(\Pi
  \begin{bmatrix}
    v_{0}^{n}\\v_{1}^{n}
  \end{bmatrix}\right)_j
  =
  \left(I - B_j^n\right)
  \begin{bmatrix}
    u^R_j\\\omega u^I_j
  \end{bmatrix}
  +
  B_j^n
  \begin{bmatrix}
    v_{0,j}^{0}\\v_{1,j}^{0}
  \end{bmatrix}
\end{align*}
so that 
  \begin{align}\label{eqn::FixedPointForm}
  \begin{bmatrix}
    v_{0,j}^{n+1} -u^R_j \\v_{1,j}^{n+1} - \omega u^I_j
  \end{bmatrix}
  =
  B_j^n
  \begin{bmatrix}
    v_{0,j}^{0} -u^R_j \\v_{1,j}^{0} - \omega u^I_j.
  \end{bmatrix}.
  \end{align}
We thus require that $B_j^n \rightarrow 0$ uniformly in $j$ to ensure convergence 
of the fixed-point iteration to the solution, $[u^R,\omega u^I]^T$, which
is true if and only if the spectral radius of $B_j$ is 
less than unity uniformly in $j$. That is, we require that 
$|\mu_j|<1$ uniformly in $j$. Defining the 
filter function $\mu(\lambda) := \beta(\lambda) + i \gamma(\lambda)$,
we may show (with a proof in Appendix~\ref{sec:FilterBounds}) 
the following lemma 
\begin{lemma}\label{lemma::FilterBounds}
  The complex-valued filter function $\mu$ satisfies $\mu(\omega) = 1$ and
  \begin{alignat*}{2}
     0 \leq    |\mu(\lambda)| &\leq
     1 - \frac{15}{32} \left(\frac{\lambda-\omega}{\omega}\right)^2, 
     &\quad{\rm when}\  \left|\frac{\lambda-\omega}{\omega}\right|\leq \frac12,\\
     |\mu(\lambda)| &\leq \frac{7}{3\pi} \approx 0.74,
     &{\rm when}\  \left|\frac{\lambda-\omega}{\omega}\right|\geq \frac12, \\
     |\mu(\lambda)| &\leq b_0\frac{\omega}{\lambda - \omega},
     &{\rm when}\  \lambda > \omega,
  \end{alignat*}
where $b_0 = 3/2\pi$. Moreover, close to $\omega$ we have the local expansion
    
  \begin{gather}\label{eqn::locexp}
    |\mu(\omega+r)| = 1- b_1 \left(\frac{r}{\omega}\right)^2 +
    R(r/\omega)\left(\frac{r}{\omega}\right)^3 ,\\ 
    b_1 = \frac{\pi^2}{6}-\frac14\approx 1.39,\quad
    ||R||_\infty\leq \frac{25 \pi^4}{4}\left(36 + 20 \pi + 250 \pi^2 + 75 \pi^3\right) \nonumber.
  \end{gather}
\end{lemma}

We denote
  \begin{align*}
    \delta_j = \frac{\lambda_j-\omega}{\omega},
  \end{align*}
the relative size of the gap between $\lambda_j$ and the Helmholtz frequency, and 
then denote the smallest gap (in magnitude) by $\delta$,
  \begin{align*}
    \delta = \delta_{j^*}, \quad j^* = {\rm argmin}_j |\delta_j|.
  \end{align*}
Then we have the following lemma
\begin{lemma}\label{lemma::SpectralRadius}
  Suppose $\delta > 0$. Then, the spectral radius $\rho$ of $\mathcal{S}$ is strictly less than one,
  and for small $\delta$,
  \begin{align}\label{eqn::rhoest}
    \rho = 1 - b_1 \delta^2 + \mathcal{O}(\delta^3),
  \end{align}
  with $b_1$ as in Lemma \ref{lemma::FilterBounds}.
  Moreover, $\mathcal{S}$ is a bounded linear map from 
  $L^2(\Omega)\times L^2(\Omega)$ to $H^1(\Omega)\times L^2(\Omega)$,
  and from $H^1(\Omega)\times L^2(\Omega)$ to $H^1(\Omega)\times H^1(\Omega)$.
\end{lemma}

\begin{proof}
Let $a_0=15/32$, $a_1=7/3\pi$ and $a_2=\sqrt{(1-a_1)/a_0}/2\approx 0.74$.
  From Lemma \ref{lemma::FilterBounds} we get
    \begin{align*}
     \rho = \sup_j |\mu(\lambda_j)|
     \leq \sup_j\max\left(1-a_0\delta_j^2,\ a_1\right)
     \leq \max\left(1-a_0\delta^2,\ a_1\right)<1.
    \end{align*}
For the more precise estimate when $\delta$ is small
we will use \eqref{eqn::locexp}. 
Since $1>\rho\geq|\mu(\omega+\omega\delta)|\to 1$ as $\delta\to 0$,
we can assume that
$\rho> 1-a_0\eta^2$, with $\eta:=\min(b_1/2||R||_\infty,\ a_2)$,
for small enough $\delta$.
Then, for $|\delta_j|>\eta$ 
we have $|\mu(\omega+\omega\delta_j)|\leq 
\max(1-a_0\eta^2/2,\ a_1)=1-a_0\eta^2/2$.
Consequently,
by Lemma \ref{lemma::FilterBounds}, we have
$$
   \rho    =
   \max_{|\delta_j|\leq \eta} |\mu(\omega+\omega\delta_j)| = 
   |\mu(\omega+\omega\delta_{k^*})|,
$$ 
for some $k^*$ with $|\delta_{k^*}|\leq \eta$. 
If $\delta_{k^*}=\delta_{j^*}$ (where $\delta=|\delta_{j^*}|$)
then \eqref{eqn::locexp} gives \eqref{eqn::rhoest}.
If not, we have
$\eta\geq |\delta_{k^*}|\geq\delta$ and by Lemma \ref{lemma::FilterBounds}
\begin{align*}
  0\leq |\mu(\omega+\omega\delta_{k^*})|-
  |\mu(\omega+\omega\delta_{j^*})|
  &=-b_1(\delta^2_{k^*}-\delta^2)
  +R(\delta_{k^*})\delta^3_{k^*}
  -R(\delta_{j^*})\delta^3_{j^*}
  \\
  &\leq
  -b_1(\delta^2_{k^*}-\delta^2)+
   \frac{b_1}2(\delta^2_{k^*}+\delta^2),
\end{align*}
which implies that $\delta_{k^*}^2\leq 3\delta^2$ and 
$\delta^2-\delta_{k^*}^2=\mathcal{O}(\delta^3)$.
Therefore
\begin{align*}
{\rho} = 
1-b_1\delta_{k^*}^2 + \mathcal{O}(\delta_{k^*}^3)
&=
1-b_1\delta^2 + 
b_1(\delta^2-\delta_{k^*}^2)+
\mathcal{O}(\delta_{k^*}^3)
\\
&=
1-b_1\delta^2 + 
\mathcal{O}(\delta_{k^*}^3+\delta^3)
\\
&=
1-b_1\delta^2 + 
\mathcal{O}(\delta^3),
\end{align*}
from which \eqref{eqn::rhoest} follows.

Letting $D:=\omega \min(1,b_0(1+1/|\delta|))$, we note that by
Lemma \ref{lemma::FilterBounds},
\begin{align*}
  |\lambda_j\mu(\lambda_j)| \leq \omega \leq D, \quad &\lambda_j\leq \omega, \\
  |\lambda_j\mu(\lambda_j)| \leq 
  \omega \frac{b_0\lambda_j}{\lambda_j-\omega} = \omega b_0(1+1/\delta_j) \leq D, \quad &\lambda_j > \omega.
\end{align*}
Moreover, the triangle inequality gives that
$|\beta(\lambda_j)|,|\gamma(\lambda_j)|  \le |\mu(\lambda_j)|$, which
implies both $\lambda_j|\beta(\lambda_j)|\le D$ and $\lambda_j|\gamma(\lambda_j)| \le D$.

Suppose now that $g,\,h \in L^2(\Omega)$ and 
  \begin{align*}
    g(x) = \sum_{j=0}^\infty g_{j} \phi_j(x), \quad h(x) = \sum_{j=0}^\infty h_{j} \phi_j(x).
  \end{align*}
Let 
$z = [g,h]^T$ and define
$C := \max \{D,|\gamma(\lambda_j)|/\lambda_j,\pi/2\omega\}$, which is bounded
for $\lambda_j \ne 0$ via the estimate \eqref{eqn::gammaBound}.
Then
 straightforward
algebra gives the bound 
\begin{align*}
   ||{\mathcal S}z||_{L^2(\Omega) \times L^2(\Omega)}^2
   &= \sum_{j=0}^\infty \left\|B_j
   \begin{bmatrix} g_j \\  h_j  \end{bmatrix}
   \right\|^2
   \\
   &= \sum_{j=0}^\infty \left|\beta(\lambda_j)g_j + \frac{\gamma(\lambda_j)}{\lambda_j} h_j\right|^2 
   +
   |\lambda_j\gamma(\lambda_j) g_j - \beta(\lambda_j)h_j|^2\\
   & 
   \leq 
   \sum_{j=0}^\infty (1 + C^2)(|g_j|^2 + |h_j|^2) 
   +
   4C |g_j| |h_j| 
   \\
   &\leq 
   \left(1 + C^2 + 2C\right) ||z||_{L^2(\Omega)\times L^2(\Omega)}^2,
\end{align*}
since $2|ab| \le a^2 + b^2$.
In the case of Neumann boundary conditions we have $\lambda_0=0$ and 
and the zeroth term must be treated specially. Using $B_0^{\rm Neu}$ we get
the same estimate
\begin{align*}
   \left\|B_0^{\rm Neu}
   \begin{bmatrix} g_0 \\  h_0  \end{bmatrix}
   \right\|^2 
   &= 
   \left|-\frac12 g_0 - \frac{\pi}{2\omega} h_0\right|^2 + \frac{1}{4} |h_0|^2
   \\
   &\le 
   \frac14 |g_0|^2 + \frac{\pi}{2\omega}|g_0||h_0| + \frac{\pi^2}{4\omega^2}|h_0|^2 + \frac14 |h_0|^2\\
   & \le 
   \left(\frac14+\frac12C\right) (|g_0|^2 + |h_0|^2) + C^2|h_0|^2
  \\
  &\le 
   \left(1 + 2C + C^2\right)(|g_0|^2 + |h_0|^2),
\end{align*}
so that in conclusion
$||{\mathcal S}z||_{L^2(\Omega) \times L^2(\Omega)}^2 \le (1 + C^2 + 2C)
||z||_{L^2(\Omega) \times L^2(\Omega)}^2$ for both cases.

Next we let
$\mathcal{S}z = [\bar{g},\bar{h}]^T$.
Then, if  $g,h \in L^2(\Omega)$, 
  \begin{align*}
||\nabla \bar{g}||_{L^2(\Omega)}^2 &\leq
    \sum_{j=0}^\infty \frac{\lambda_j^2}{c^2_\text{min}}\left|\beta(\lambda_j)g_j + \frac{\gamma(\lambda_j)}{\lambda_j} h_j\right|^2 
    \le
    \sum_{j=0}^\infty \frac{2}{c^2_\text{min}}\left(|\lambda_j\beta(\lambda_j)g_j|^2 + 
    |\gamma(\lambda_j)h_j|^2\right)
\\
    &\le 
      \sum_{j=0}^\infty \frac{2(D^2+1)}{c^2_\text{min}}\left(|g_j|^2 
          + |h_j|^2 \right) 
    = \frac{2(D^2 + 1)}{c^2_\text{min}} \|z\|^2_{L^2(\Omega) \times L^2(\Omega)},
  \end{align*}
which gives
  \begin{align*}
    \|\mathcal{S}z\|_{H^1(\Omega) \times L^2(\Omega)}^2 
    &= 
    \|\mathcal{S}z\|_{L^2(\Omega) \times L^2(\Omega)}^2
    +||\nabla \bar{g}||_{L^2(\Omega)}^2
    \\
    &\leq
  \left( 1 +C^2 + 2C + \frac{2(D^2 + 1)}{c^2_\text{min}} \right)
    \|z\|^2_{L^2(\Omega) \times L^2(\Omega)},
  \end{align*}
showing that $\mathcal{S}$ is a bounded linear map 
from $L^2(\Omega)\times L^2(\Omega)$ to $H^1(\Omega)\times L^2(\Omega)$.

If instead $g \in H^1(\Omega)$ and $h \in L^2(\Omega)$, we compute
  \begin{align*}
||\nabla \bar{h}||_{L^2(\Omega)}^2 
    &\leq
    \sum_{j=0}^\infty \frac{\lambda_j^2}{c^2_\text{min}}|-\gamma(\lambda_j)\lambda_jg_j + \beta(\lambda_j)h_j|^2 
    \\
    &\le 
        \sum_{j=0}^\infty \frac{2}{c^2_\text{min}}\left(
        |\gamma(\lambda_j)\lambda_j^2g_j|^2 +|\beta(\lambda_j)\lambda_jh_j|^2\right)\\
    &\le 
        \frac{2D^2}{c^2_\text{min}}
        \sum_{j=0}^\infty \left(\lambda_j^2|g_j|^2 
          + |h_j|^2 \right)
          \\
          &\leq \frac{2D^2}{c^2_\text{min}}
          \left(c^2_\text{max}\|\nabla g\|^2_{L^2(\Omega)}+ \|h\|_{L^2(\Omega)}\right)
          \leq
          \frac{2D^2(c^2_\text{max}+1)}{c^2_\text{min}}
    ||z||_{H^1(\Omega) \times L^2(\Omega)}^2           
  \end{align*}
so that
  \begin{align*}
      \|\mathcal{S}z\|_{H^1(\Omega) \times H^1(\Omega)}^2 
    &= 
    \|\mathcal{S}z\|_{H^1(\Omega) \times L^2(\Omega)}^2
    +||\nabla \bar{h}||_{L^2(\Omega)}^2
    \\
    &\leq     
    \left( 1 +C^2 + 2C + \frac{2D^2(c^2_\text{max}+2)+2}{c^2_\text{min}} 
          \right)
    \|z\|^2_{H^1(\Omega) \times L^2(\Omega)}
  \end{align*}
which shows that $\mathcal{S}$ is a bounded linear map from $H^1(\Omega)\times L^2(\Omega)$
to $H^1(\Omega) \times H^1(\Omega)$, concluding the proof of the lemma. 
\end{proof}

Further, denoting $e^n := [\text{Re}\{u\} - v_0^n, \, \omega\text{Im}\{u\} - v_1^n]^T = [e_0^n, e_1^n]^T$, 
from 
\eqref{eqn:iteration} 
we obtain
  \begin{align*}
    e^{n} = \mathcal{S}[\text{Re}\{u\} - v_0^{n-1},\, \omega\text{Im}\{u\} - v_1^{n-1}]^T
            = \mathcal{S} e^{n-1}
            = \mathcal{S}^n e^0,
  \end{align*}
which shows that $e^n \rightarrow 0$ since ${\mathcal S}^n \rightarrow 0$. Thus the
iterates $[v_0^n, \, v_1^n]^T$ converge to $[\text{Re}\{u\}, \, \omega\text{Im}\{u\}]^T$
in $L^2(\Omega)\times L^2(\Omega)$. Since $v_0^0= v_1^0 = 0$, it follows 
from Lemma \ref{lemma::SpectralRadius} 
that the iterates
$[v_0^n, \, v_1^n]^T$ and $[e_0^n, \, e_1^n]^T$ both belong to
$H^1(\Omega) \times H^1(\Omega)$ for $n>0$
when $u\in H^1(\Omega)$, 
as $\mathcal{S}$ is a bounded linear map from $H^1(\Omega)\times L^2(\Omega)$
to $H^1(\Omega) \times H^1(\Omega)$.
We can therefore also get convergence in $H^1(\Omega) \times H^1(\Omega)$.
To show this, let
  \begin{align*}
    \beta_j + i\gamma_j = r_j \exp (i \phi_j), \quad r_j^2 = |\beta_j|^2 + |\gamma_j|^2, \quad \phi_j = \arctan(\gamma_j/\beta_j).
  \end{align*}
It can then be shown that powers of the operator $B_j$ can be written as 
\begin{align*}
   B_j^n =  
   r_j^n
   \begin{pmatrix}
    \cos(n\phi_j) 
    & \sin(n\phi_j)/\lambda_j  \\
    -\lambda_j \sin(n\phi_j) 
    & \cos(n \phi_j)
   \end{pmatrix},
\end{align*}
where each entry is bounded and goes to zero in the limit as $n\rightarrow \infty$
since the spectral radius of $B_j$ is less than one. 
Then, using the fact that $|r_j\lambda_j|\leq D$,
  \begin{align*}
  ||\nabla e^n||^2_{L^2(\Omega)\times L^2(\Omega)} &=
    ||\nabla \mathcal{S}^n  e^0||^2_{L^2(\Omega)\times L^2(\Omega)} \\
    &
    \le
    \sum_{j=0}^\infty \frac{\lambda_j^2r_j^{2n}}{c^2_\text{min}}\left(\left|\cos(n\phi_j)e_{j,0}^0 +\frac{\sin(n\phi_j)}{\lambda_j} e_{j,1}^1\right|^2 \right. 
    \\ 
    &\left.\quad\quad\quad\quad\quad\quad\quad+
    \left|-\lambda_j \sin(n\phi_j)e_{j,0}^0 + \cos(n\phi_j) e_{j,1}^1\right|^2 \right)\\
    &\le
    \sum_{j=0}^\infty \frac{2r_j^{2n}}{c^2_\text{min}}\left(
    \lambda_j^2 |e_{j,0}^0|^2 + |e_{j,1}^1|^2 \right)
    +  \frac{2r_j^{2n}\lambda_j^2}{c^2_\text{min}}\left(\lambda_j^2|e_{j,0}^0|^2 
                                                            +|e_{j,1}^1|^2 \right)  \\
    &\le \left(\sup_{j} r_j\right)^{2n-2}
    \frac{2(1 + D^2)}{c^2_\text{min}}
\sum_{j=0}^\infty \left(\lambda_j^2 |e_{j,0}^0|^2 
 + |e_{j,1}^1|^2 \right)\\
    & \le  \rho^{2n-2} \frac{2(1 + D^2)}{c^2_\text{min}}
\left(c^2_\text{max} \|\nabla e_0^0\|_{L^2(\Omega)}^2 
    +\|e_1^0\|_{L^2(\Omega)}^2 \right) \rightarrow 0.
  \end{align*}
We conclude that the iteration converges in $ H^1(\Omega) \times H^1(\Omega)$
with convergence rate $\rho$. By Lemma \ref{lemma::SpectralRadius} we 
have $\rho \sim 1 - 1.39 \delta^2$ so that the smallest gap, $\delta$,
determines the convergence rate. We thus have proven the following theorem
\begin{theorem}\label{thm::FirstOrderConv}
Suppose $u\in H^1(\Omega)$ is the solution
of the Helmholtz equation 
(\ref{eq:helm}).
The iteration in \eqref{eqn::Iteration} and \eqref{eq::filterstep} converges
in $H^1(\Omega) \times H^1(\Omega)$
for the Dirichlet and Neumann problems
away from resonances
to 
$[\text{Re}\{u\}, \, \omega\text{Im}\{u\} ]^T$.
The convergence rate is $1-\mathcal{O}(\delta^2)$,
where $\delta$ is the minimum gap between $\omega$ and 
the eigenvalues of $-\nabla\cdot(c^2(x)\nabla)$.
\end{theorem}


\subsection{Convergence in the Non-Energy Conserving Case}\label{sec::non-energy}


With Theorem~\ref{thm::FirstOrderConv} providing convergence of the general 
WaveHoltz iteration in the energy conserving case, we turn toward proving 
convergence for problems with impedance boundary conditions.
For simplicity we prove convergence in a single spatial 
dimension. We note that it is possible to use
the following approach to prove convergence for certain problems in 
higher dimensions, e.g. problems with a constant wavespeed
in certain simple geometries.
Consider now the following Helmholtz problem with impedance boundary conditions
  \begin{align}\label{eqn::HelmImp1D}
    \frac{\partial}{\partial x}\left[c^2(x) \frac{\partial}{\partial x}u(x)\right] & + \omega^2 u(x) = f(x), \quad a \le x \le b, \nonumber\\
    i \alpha \omega u(a) & - \beta c(x) u_x(a) = 0, \\ 
    i \alpha \omega u(b) & + \beta c(x) u_x(b) = 0\nonumber,
  \end{align}
where $\alpha^2+\beta^2=1$, but
$\alpha, \, \beta \ne 0$. 
As before we assume that
$c\in L^\infty(a,b)$ with
the bounds $0< c_{\min}\leq c(x)\leq c_{\max}<\infty$ a.e. in $(a,b)$,
and $f \in L^2(a,b)$, but we now
require additionally that $f$ is compactly supported in $(a,b)$
and
that $c(a) = c_a$, $c(b) = c_b$ with $c$ constant in a neighborhood of the endpoints. 
We reformulate this in the time domain as 
  \begin{eqnarray*}
    &&  w_{tt} =   \frac{\partial}{\partial x} \left[c^2(x) \frac{\partial}{\partial x}w\right]
      - f(x) e^{-i \omega t}, \quad a\le x \le b, \ \ 0 \le t \le T, \nonumber \\
    &&  w(0,x) = v_0(x), \quad w_t(0,x) = v_1(x), \nonumber \\
    &&  \alpha w_t(t,a) -\beta c(x) w_x(t,a)=0,\\
    &&  \alpha w_t(t,b) +\beta c(x) w_x(t,b)=0.
  \end{eqnarray*}
In general, the solution of the above equation will yield complex-valued solutions 
and so we take the real part of the 
equation as shown earlier and use the general iteration \eqref{eqn::Iteration}. 
Note that in $1$D  the impedance boundary conditions with $\alpha = \beta = 1/\sqrt{2}$ are 
equivalent to outflow/radiation conditions when the initial data is compactly supported
in the interval $[a,b]$.
If $\alpha \ne \beta$ with $\alpha, \, \beta \ne 0$, then in addition to outgoing waves at the boundary
there will be reflections due to the impedance boundary condition.
In either case, if we let $\tilde a < a - c_aT/2$ and $\tilde b > b + c_bT/2$, 
then $w$ is equal to $\tilde w$ on $[a,b]$ for $t \in [0,T]$
if $\tilde w$ solves the following Neumann problem
in the extended domain $\widetilde \Omega := [\tilde a, \tilde b] \supset [a,b] := \Omega$
  \begin{eqnarray}\label{eqn::ExtTD}
    &&  \tilde w_{tt} =   \frac{\partial}{\partial x} \left[\tilde{c}^2(x) \frac{\partial}{\partial x}\tilde w\right]
      - \text{Re}\{\tilde f(x) e^{-i\omega t}\}, \quad \tilde a\le x \le \tilde b, \ \ 0 \le t \le T, \nonumber \\
    &&  \tilde w(0,x) = \tilde{v}_0(x), \quad \tilde w_t(0,x) = \tilde{v}_1(x), \\
    &&  \tilde w_x(t,\tilde a) = 0, \quad \tilde w_x(t,\tilde b)=0 \nonumber,
  \end{eqnarray}
where $\tilde{v}_0$ and $\tilde c$ are the constant extensions
(with $\gamma = \alpha/\beta$)
  \begin{align*}
    \tilde{v}_0(x) = 
      \begin{cases}
        v_0(a), & \tilde a \le x < a, \\
        v_0(x), & a \le x \le b, \\
        v_0(b), & b < x \le \tilde b,
      \end{cases}
      \quad 
      \tilde{c}(x) = 
        \begin{cases}
          \gamma c_a, & \tilde a \le x < a, \\
          c(x), & a \le x \le b, \\
          \gamma c_b, & b < x \le \tilde b,
        \end{cases}
    \end{align*}
and $\tilde{v}_1, \tilde f$ are zero extensions of $v_1$ and $f$,
  \begin{align*}
    \tilde{v}_1(x) = 
      \begin{cases}
        0, & \tilde a \le x < a, \\
        v_1(x), & a \le x \le b, \\
        0, & b < x \le \tilde b,
      \end{cases}
      \quad 
      \tilde{f}(x) = 
        \begin{cases}
          0, & \tilde a \le x < a, \\
          f(x), & a \le x \le b, \\
          0, & b < x \le \tilde b.
        \end{cases}
    \end{align*}
%

That is, we extend the domain such that traveling waves may reflect 
off of the Neumann boundary but not re-enter the domain of interest,
$a \le x \le b$, within a period $T$ 
(see Appendix~\ref{sec:WaveExtension} for an outline of the construction).
Let $\Pi$ be the WaveHoltz integral operator
\eqref{eq::filterstep} on the original domain $\Omega$
with impedance boundary conditions.
We recall that iterates generated by 
$\Pi$ at a given point, $x \in \Omega$, are 
the time-average of the wave solution at $x$ generated by
the input data.
Since the extended wave solution $\tilde w(t,x) = w(t,x)$ 
for $0 \le t \le T$, we may write $\Pi = P \widetilde{\Pi} E$
where $P$ is a projection operator onto the initial interval, i.e.
$P v(x) = v(x)|_{a\le x \le b}$, $E$ is the 
extension operator such that $[v_0,v_1]^T \rightarrow [\tilde v_0,\tilde v_1]^T$,
and $\widetilde{\Pi}$ is the WaveHoltz operator on the domain $\widetilde{\Omega}$.
If it can be guaranteed that $\omega^2 \ne \lambda_j^2$ where 
$\lambda_j^2$ is an eigenvalue of the operator $-\partial_x(\tilde c^2(x) \partial_x)$,
then we may prove convergence as was done for Theorem~\ref{thm::FirstOrderConv}.

To show this, results on the continuity of eigenvalues of the Laplacian from \cite{kong1996eigenvalues}
will be used. We present the framework of \cite{kong1996eigenvalues} 
needed here and consider the following differential equation
  \begin{align}\label{eqn::SL_eq}
    -\frac{d}{dx}\left(c^2 \frac{d}{dx}y(x)\right) = \lambda y(x), \quad x \in (a',b'), \quad -\infty \leq a' < b' \le \infty, \quad \lambda \in \mathbb{R},
  \end{align}
where $c^2:(a',b') \rightarrow \mathbb{R}$ and  $1/c^2 \in L^1_{\text{loc}}(a',b')$.
Letting $I=[a,b], \quad a'<a<b<b'$ and additionally imposing the Neumann 
conditions $y'(a) = 0 = y'(b)$, the above Sturm-Liouville (SL) problem
is such that all eigenvalues are real, simple, and can be ordered to satisfy
  \begin{align}
     0 \le \lambda^2_0 < \lambda^2_1 < \lambda^2_2 < \dots; \quad \lim_{n\rightarrow \infty} \lambda^2_n = + \infty.\label{eqn::LamOrder}
  \end{align}

Under the above assumptions, we state the following theorem that is proven in \cite{kong1996eigenvalues}.
\begin{theorem}[Kong \& Zettle]\label{thm::EvalDiff}
  Let  $1/c^2 \in L^1_{\text{loc}}(a',b')$, fix $a',b'$, 
  and suppose $a,b$ are such that $a'<a<b<b'$.
  Fix $a$ and let $\lambda_n(b)$ be an eigenvalue of the 
  SL problem \eqref{eqn::SL_eq} with homogeneous Neumann 
  boundary conditions at $x=a$ and $x=b$ with the
  corresponding eigenfunction $u_n(x;b)$.
  Then the eigenvalue $\lambda_n \in C^1([a',b'])$ satisfies 
  the following differential equation:
    \begin{align*}
      \frac{d}{db}\lambda_n(b) = -\lambda_n(b)u_n^2(b;b).
    \end{align*}
\end{theorem}

That is, the eigenvalues of the SL problem \eqref{eqn::SL_eq} are differentiable
functions of the endpoint $b$. This gives us the following useful corollary.
  \begin{corollary}{}{}\label{corr::EigDiff}
    For $n=1,2,\dots$, $\lambda_n(b)$ is a strictly decreasing function of 
    $b$ on $[a',b']$.
  \end{corollary}
  \begin{proof}
    For homogeneous Neumann conditions, we have
    that $u_n'(b;b) = 0$. It follows that $u_n(b;b) \ne 0$ as otherwise 
    $u_n(b;b) \equiv 0$ since $u_n$ satisfies a linear, homogeneous 
    second order ODE. As $\lambda_n(b)>0$ for $n>0$ we then have
      \begin{align*}
        \frac{d}{db}\lambda_n(b) = -\lambda_n(b) u_n^2(b;b) < 0,
      \end{align*}
    so that $\lambda_n(b)$ is a strictly decreasing function of the endpoint $b$.
  \end{proof}

As a consequence of Theorem~\ref{thm::EvalDiff}, we have
\begin{lemma}\label{lemma:Impedance}
  Suppose $\omega >0$ and that we extend $c$ to $(a',b')$ with
  $a'<a-c_aT/2$ and $b'>b + c_bT/2$.
Fix $\tilde a \in (a',a - c_aT/2)$.
Then there exists an endpoint $\tilde b \in (b + c_bT/2,b')$ such that $\omega^2 \ne \lambda_n(\tilde b)$ for each $n \in \mathbb{N}_0$, where $\lambda_n$
are the Neumann eigenvalues.
\end{lemma}
\begin{proof}
We note first that $\tilde{c}\in L^\infty(a',b')$
with $\tilde{c}(x)\geq \min(1,\gamma) c_\text{min}$. Hence, $1/\tilde{c}^2\in L^1_{loc}(a',b')$, so 
Theorem~\ref{thm::EvalDiff} and
Corollary~\ref{corr::EigDiff} apply.
  Clearly we have $\lambda_0(t) = 0$ for every $t$, and since $\omega > 0$ we have
  $\omega^2 \ne \lambda_0(t)$. Suppose now that $t\in(b+c_bT/2,b')$ is such 
  that $\omega^2 = \lambda_n(t)$ 
  for some $n\in \mathbb{N}$. (If not, we take $\tilde{b}=t$.)
  Recall that by \eqref{eqn::LamOrder} we have that $\omega^2 = \lambda_n(t) < \lambda_{n+1}(t)$. 
  Since $\lambda_n(t), \lambda_{n+1}(t)$ are continuous, decreasing functions 
  of the endpoint by Corollary~\ref{corr::EigDiff}, 
  there necessarily exists $\delta \in(0,b'-t)$ such that
    \begin{align*}
      \lambda_n( t+\delta) < \omega^2 < \lambda_{n+1}( t+\delta).
    \end{align*}
  Letting $\tilde b = t + \delta$ we thus have that $\omega^2 \ne \lambda_n(\tilde b)$ for each $n \in \mathbb{N}_0$, as desired.
\end{proof}

From this we can prove the following theorem, in which 
we demonstrate convergence in \mbox{$H^1(\Omega) \times L^2(\Omega)$}
rather than \mbox{$H^1(\Omega) \times H^1(\Omega)$}. 
\begin{theorem}\label{thm:TheoremOutflow}
  Let the 1D domain $\Omega = [a,b]$ be a bounded interval.
  Suppose $f\in L^2(\Omega)$ is compactly supported in $\Omega$, and 
  $c\in L^\infty(a,b)$ with
the bounds $0< c_{\min}\leq c(x)\leq c_{\max}<\infty$ a.e. in $(a,b)$,
and the additional restriction that
%
%
  $c(a) = c_a$, $c(b) = c_b$, with $c$ constant near the endpoints.
    Under these conditions, if $u\in H^1(\Omega)$ is the solution
of the Helmholtz problem with impedance boundary conditions \eqref{eqn::HelmImp1D},
 the iteration \eqref{eqn::Iteration} and \eqref{eq::filterstep} converges in $H^1(\Omega) \times L^2(\Omega)$ to 
$[\text{Re}\{u\}, \, \omega\text{Im}\{u\} ]^T$. 
\end{theorem}
\begin{proof}

  By Lemma \ref{lemma:Impedance}, there exists an extended wave equation \eqref{eqn::ExtTD} on the domain
  $\widetilde{\Omega} = [\tilde a, \tilde b]$
  with homogeneous Neumann boundary conditions such that the eigenvalues $\tilde \lambda_j$ of the
  Laplacian, $-\partial_x(\tilde{c}^2 \partial_x)$, on $\widetilde{\Omega}$ are not in resonance. 
  Defining $\tilde \beta_j = \beta(\tilde \lambda_j)$, $\tilde \gamma_j = \gamma(\tilde \lambda_j)$,
  and $\tilde \mu_j = \mu(\tilde \lambda_j)$, this immediately gives that the 
  spectral radius of the WaveHoltz operator, $\tilde \rho = \sup_j |\tilde \mu_j|$, is smaller than one.
  Moreover, the extended wave solution $\tilde w$ on $\widetilde{\Omega}$ coincides with the interior impedance wave 
  solution $w$ on $\Omega$ for $t \in [0,T]$.
  For the extended speed function $\tilde{c}$ we have the bounds
  $0< \tilde{c}_{\min}\leq \tilde{c}(x)\leq \tilde{c}_{\max}<\infty$ a.e. in 
$(\tilde{a},\tilde{b})$, where $\tilde c_\text{min}=\min(1,\gamma) c_\text{min}$
and
    $\tilde c_\text{max}=\max(1,\gamma) c_\text{max}$.

  Letting $u$ be the solution of the Helmholtz equation \eqref{eqn::HelmImp1D}, 
  we define $q(t,x) = \cos(\omega t)[\text{Re}\{u\}, \, \omega\text{Im}\{u\}]^T$ the
  time-harmonic Helmholtz solution in $\Omega$ and $\tilde w^n(t,x)$ the solution of 
  \eqref{eqn::ExtTD} with initial data $\tilde v_0^n, \tilde v_1^n$. Letting the error be 
  $e^n := [\text{Re}\{u\} - v_0^n, \, \omega\text{Im}\{u\} - v_1^n]^T = [e_0^n, e_1^n]^T$,
  it is clear the difference \mbox{$d(t,x) = q(t,x) - w(t,x)$} satisfies the unforced,
  homogeneous wave equation
  \begin{eqnarray*}
    &&   d_{tt} =   \frac{\partial}{\partial x} \left[{c}^2(x) \frac{\partial}{\partial x} d\right]
      , \quad  a\le x \le  b, \ \ 0 \le t \le T, \nonumber \\
    &&   d(0,x) = {e}_0(x), \quad  d_t(0,x) = {e}_1(x), \\
    &&  \alpha d_t(t,a) -\beta c(x) d_x(t,a)=0,\\
    &&  \alpha d_t(t,b) +\beta c(x) d_x(t,b)=0.
  \end{eqnarray*}
  It follows that the WaveHoltz iteration applied to the error is of the form
  \begin{equation}\label{eqn::ExtendedIterationSpectral}
    e^{n+1}
     = \Pi
    e^n=
    P\widetilde{\mathcal{S}}E
e^n    =
    (P\widetilde{\mathcal{S}}E)^{n+1}
e^0,
  \end{equation}
  where $\widetilde{\mathcal{S}}$ is 
  defined in \eqref{eqn::SpectralOp}, but
  with respect to the eigenbasis of the extended Laplacian.
  Note that $e^0 =[\text{Re}\{u\}, \, \omega\text{Im}\{u\}]^T
  \in H^1(\Omega)\times H^1(\Omega)$ since
  $v_0^0 = v_1^0 = 0$.
    
We note further that the extension operator $E$ maps
$H^1(\Omega)\times L^2(\Omega)$
to $H^1(\tilde\Omega)\times L^2(\tilde\Omega)$, while
for the projection operator we have $P:H^1(\tilde\Omega)\times H^1(\tilde\Omega)\mapsto H^1(\Omega)\times L^2(\Omega)$
and the bound
  \begin{align}\label{eqn:Pbound}
     ||Pz||_{H^1(\Omega)\times L^2(\Omega)}\leq
     ||z||_{H^1(\tilde\Omega)\times L^2(\tilde\Omega)}.
  \end{align}
    In Lemma~\ref{lemma::SpectralRadius} it 
    was shown that 
    $\widetilde{\mathcal{S}}:H^1(\tilde\Omega)\times L^2(\tilde\Omega)\mapsto 
    H^1(\tilde\Omega)\times H^1(\tilde\Omega)$ and it
    follows that
  \begin{align}\label{eqn:Smaps}
    \widetilde{\mathcal{S}}E: H^1(\Omega)\times L^2(\Omega)\mapsto 
    H^1(\tilde\Omega)\times H^1(\tilde\Omega), \\
    \widetilde{\mathcal{S}}EP:
    H^1(\tilde\Omega)\times L^2(\tilde\Omega)\mapsto H^1(\tilde\Omega)\times H^1(\tilde\Omega).\nonumber
  \end{align}
We define
  \begin{align*}
  \tilde e^{n+1} = \widetilde{\mathcal{S}}EP\tilde e^n, \qquad \tilde e^0 = \widetilde{\mathcal{S}}Ee^0, \qquad
  \tilde e^n = 
      \left[
    \begin{array}{c}
    \tilde e_0^{n}\\
    \tilde e_1^{n} 
    \end{array}
    \right].
  \end{align*}
Since $e^0 \in H^1(\Omega) \times H^1(\Omega)$ it follows 
from \eqref{eqn:Smaps} that
$\tilde e^n \in H^1(\tilde\Omega) \times H^1(\tilde\Omega)$ for all $n\geq 0$.
Moreover, by rearranging the iteration \eqref{eqn::ExtendedIterationSpectral} we obtain
  $$
    e^{n+1}=(P\widetilde{\mathcal{S}}E)^{n+1} e^0 = P (\widetilde{\mathcal{S}}EP)^{n} \widetilde{\mathcal{S}}Ee^0 = P\tilde e^{n+1}.
$$
Then by \eqref{eqn:Pbound} we have for $n\geq 1$,
$$
   ||e^n||_{H^1(\Omega)\times L^2(\Omega)}\leq
   ||\tilde e^n||_{H^1(\tilde\Omega)\times L^2(\tilde\Omega)},
$$
and to prove the stated convergence it is therefore sufficent to prove
that $\tilde e^n\to 0$ in $H^1(\tilde\Omega)\times L^2(\tilde\Omega)$.

For the convergence we consider first the energy semi-norm
  $\|\cdot\|_c$ on $H^1(\widetilde{\Omega}) \times L^2(\widetilde{\Omega})$.
Let $\tilde z = [\tilde v_0, \tilde v_1]^T\in H^1(\widetilde{\Omega}) \times L^2(\widetilde{\Omega})$ and define
\begin{align*}
      \|\tilde z\|^2_c := 
      \left\|\tilde c \frac{\partial \tilde v_0}{\partial x} \right\|_{L^2(\widetilde{\Omega})}^2 
       + \|\tilde v_1\|_{L^2(\widetilde{\Omega})}^2
       =\sum_{j=0}^\infty \tilde \lambda_j^2 | \tilde v_{0,j}|^2 +| \tilde v_{1,j}|^2.
    \end{align*}
In this semi-norm we have that 
  \begin{align*}
    \|EP\tilde z\|_c^2
          = \int_{\widetilde{\Omega}} \left|\tilde c \frac{\partial }{\partial x} EP \tilde v_0\right|^2 + |EP\tilde v_1|^2  \, dx
          &= \int_{\Omega} \left|\tilde c \frac{\partial \tilde v_0}{\partial x} \right|^2 + |\tilde v_1|^2  \, dx
          \\
          &\le \int_{\widetilde{\Omega}} \left|\tilde c \frac{\partial \tilde v_0}{\partial x} \right|^2 + |\tilde v_1|^2  \, dx
          \le \|\tilde z\|_c^2.
  \end{align*}

  We now proceed with the proof and
  define $\tilde y = EP \tilde z$, where $\tilde y$ has the form
    \begin{align*}
      \tilde y = EP\sum_{j=0}^\infty \begin{bmatrix} \tilde v_{0,j}\\ \tilde v_{1,j} \end{bmatrix}\phi_j
          = \sum_{j=0}^\infty \begin{bmatrix}\tilde y_{0,j} \\ \tilde y_{1,j} \end{bmatrix}\phi_j.
    \end{align*}
  It follows that
    \begin{align*}
      \widetilde{\mathcal{S}}EP \tilde z = \sum_{j=0}^\infty B_j\begin{bmatrix}\tilde y_{0,j} \\ \tilde y_{1,j} \end{bmatrix}\phi_j
      =
      \sum_{j=1}^\infty 
      \begin{bmatrix} 
        \tilde \beta_j \tilde y_{0,j} + \tilde \gamma_j \tilde y_{1,j}/\tilde \lambda_j\\ 
        -\tilde \lambda_j \tilde y_{0,j} + \tilde \beta_j \tilde y_{1,j}  
      \end{bmatrix} \phi_j
+ B_0^{\rm Neu}
\begin{bmatrix}\tilde y_{0,0} \\ \tilde y_{1,0} \end{bmatrix}\phi_0,
    \end{align*}
  so that
    \begin{align*}
      \|\widetilde{\mathcal{S}}EP \tilde z\|_c^2 = \sum_{j=1}^\infty \tilde \lambda_j^2\left(\tilde \beta_j \tilde y_{0,j} + \frac{\tilde \gamma_j}{\tilde \lambda_j} \tilde y_{1,j}\right)^2
      +
      \sum_{j=0}^\infty (-\tilde \lambda_j \tilde \gamma_j\tilde y_{0,j} + \tilde \beta_j \tilde y_{1,j})^2.
    \end{align*}
  Since $\tilde \beta_j^2 + \tilde \gamma_j^2 = |\tilde \mu_j|^2 \le \tilde \rho^2 < 1$,
  a simple expansion shows that
    \begin{align*}
      \|\widetilde{\mathcal{S}}EP\tilde z\|_c^2 = \sum_{j=0}^\infty (\tilde \beta_j^2 + \tilde \gamma_j^2) (\tilde \lambda_j^2 |\tilde y_{0,j}|^2 + |\tilde y_{1,j}|^2)
      &\le 
      \left(\sup_j |\tilde \mu_j|^2\right) \sum_{j=0}^\infty \tilde \lambda_j^2 |\tilde y_{0,j}|^2 + |\tilde y_{1,j}|^2
      \\ 
      &\le \tilde \rho^2 \|\tilde y\|_c^2
      \le \tilde \rho^2 \|\tilde z\|_c^2.
    \end{align*}
  We then obtain the estimate
    \begin{align*}
\|\tilde e^n\|_c^2=      \|(\widetilde{\mathcal{S}}EP)^n \tilde e^0\|_c^2 \le \tilde \rho^2  \|(\widetilde{\mathcal{S}}EP)^{n-1} \tilde e^{0}\|_c^2 \le \dots \le \tilde \rho^{2n} \|\tilde e^0\|_c^2 \rightarrow 0.
    \end{align*}
    
We now consider the full 
$H^1(\tilde\Omega)\times L^2(\tilde\Omega)$-norm.
  An application of the triangle and Poincar\'e inequality (with constant $C_p$) gives
     \begin{align}\label{eqn::Poincare}
      \|\tilde e^n\|^2_{H^1(\widetilde{\Omega})\times L^2(\widetilde{\Omega})} 
      &= \|\partial_x \tilde e_0^n\|^2_{L^2(\widetilde{\Omega})}
      +\|\tilde e_0^n\|^2_{L^2(\widetilde{\Omega})} + \|\tilde e_1^n\|^2_{L^2(\widetilde{\Omega})} 
      \\
      &\le  \|\partial_x \tilde e_0^n\|^2_{L^2(\widetilde{\Omega})} + \|\tilde e_0^n -\tilde e^n_{0,0}\phi_0\|^2_{L^2(\widetilde{\Omega})} +\|\tilde e^n_{0,0}\phi_0\|^2_{L^2(\widetilde{\Omega})} + \|\tilde e_1^n\|^2_{L^2(\widetilde{\Omega})} \nonumber \\ 
      &\le (1+C_p) \|\partial_x \tilde e_0^n\|^2_{L^2(\widetilde{\Omega})} + 
      |\tilde e^n_{0,0}|^2+ \|\tilde e_1^n\|^2_{L^2(\widetilde{\Omega})} ,
      \nonumber\\
      &\le (1+(1+C_p)\tilde c_{\rm min}^{-2}) \|\tilde e^n\|^2_{c} + |\tilde e^n_{0,0}|^2,
      \nonumber
    \end{align}
%
  where  $\phi_0$ is a constant eigenfunction
  of the Laplacian (and thus of $\widetilde{\mathcal{S}}$) with
  eigenvalue $\lambda_0 = 0$.
  To obtain convergence in $H^1(\widetilde \Omega)\times L^2(\widetilde \Omega)$
  of the error $\tilde e^n$ we must
 thus  examine the convergence of $\tilde e^n_{0,0}$ separately.

  Before proceeding, we require the following lemma:
  \begin{lemma}\label{lemma::L2toEnergyEstimate}
  Let $\tilde z \in H^1(\widetilde{\Omega}) \times L^2(\widetilde{\Omega})$. Then
  \begin{align*}
    \|EP\tilde z - \tilde z\|_{L^2(\widetilde{\Omega}) \times L^2(\widetilde{\Omega})}^2 \le C_1 \|\tilde z\|_c^2,
  \end{align*}
  where $C_1 = \max\{2(a-\tilde a)/\tilde c^2_\text{min},2(\tilde b - b)/\tilde c^2_\text{min},1\} > 0$.
  \end{lemma}
  \begin{proof}
    Let $\tilde z = [\tilde v_0, \tilde v_1]^T$ with $\tilde v_0 \in H^1(\widetilde{\Omega})$ and $\tilde v_1 \in L^2(\widetilde{\Omega})$. Then
    \begin{align*}
      \|EP\tilde z - \tilde z\|_{L^2(\widetilde{\Omega}) \times L^2(\widetilde{\Omega})}^2
      & =
      \int_{\tilde a}^a |\tilde v_0(x) - \tilde v_0(a)|^2 + |\tilde v_1(x)|^2  \, dx
      \\
      &\quad +
      \int_b^{\tilde b} |\tilde v_0(x) - \tilde v_0(b)|^2 + |\tilde v_1(x)|^2  \, dx
      \\
      & \le 
      \int_{\tilde a}^a \left|\int_a^x \partial_x\tilde v_0(s) \, ds\right|^2   \, dx
      +
      \int_b^{\tilde b} \left|\int_b^x \tilde \partial_xv_0(s) \, ds\right|^2   \, dx
      \\
      & \quad + \|\tilde v_1\|_{L^2(\widetilde{\Omega})},
    \end{align*}
  since $\tilde v_0 \in H^1(\widetilde\Omega)$. 
Moreover,
    \begin{align*}
      \int_{\tilde a}^a \left|\int_a^x \partial_x\tilde v_0(s) \, ds\right|^2   \, dx
      \le 
      \int_{\tilde a}^a \int_a^x |\partial_x\tilde v_0(s)|^2 \, ds   \, dx
      &\le
      \int_{\tilde a}^a \frac{1}{\tilde c^2_\text{min}}\|\tilde c\partial_x\tilde v_0\|_{L^2(\widetilde\Omega)}^2 \, dx
      \\
      &\le 
      \frac{a-\tilde a}{\tilde c^2_\text{min}}\|\tilde c\partial_x\tilde v_0\|_{L^2(\widetilde\Omega)}^2.
    \end{align*}
  A similar estimate for the integral in the left part of the extended domain, $x \in [b,\tilde b]$,
  gives the bound
    \begin{align*}
      \|EP\tilde z - \tilde z\|_{L^2(\widetilde{\Omega}) \times L^2(\widetilde{\Omega})}^2
      \le 
      C_1\|\tilde c\partial_x\tilde v_0\|_{L^2(\widetilde\Omega)}^2
      + 
      \|\tilde v_1\|_{L^2(\widetilde\Omega)}^2
      \le 
      C_1 \|\tilde z\|_c^2,
    \end{align*}
  as desired.
  \end{proof}

To simplify notation we define
    \begin{align*}
      I_0      \begin{bmatrix}
        \tilde v_0 \\ \tilde v_1
      \end{bmatrix}
 := 
      \tilde v_0,\qquad
      I^*_0  \tilde v_0:=    \begin{bmatrix}
        \tilde v_0 \\ 0
      \end{bmatrix},\qquad
      \langle f,g\rangle :=\int_{\tilde{a}}^{\tilde{b}}f(x)g(x)dx,
    \end{align*} 
  so that we may write the constant component of the $\tilde{e}^n_0$ error as 
    \begin{align*}
       \tilde e_{0,0}^{n+1} 
       = \langle I_0 \tilde e^{n+1}, \phi_0\rangle 
       = \langle I_0 \widetilde{\mathcal{S}} E P \tilde e^n, \phi_0\rangle 
        &=
       \langle I_0 \widetilde{\mathcal{S}} (E P \tilde e^n - \tilde e^n), \phi_0\rangle 
            \\
            & + \langle I_0 \widetilde{\mathcal{S}}(I-I_0^*I_0) \tilde e^n, \phi_0\rangle
            + \langle I_0 \widetilde{\mathcal{S}}I_0^*I_0 \tilde e^n, \phi_0\rangle.
     \end{align*}
     For the last term we get, since $B_0=B_0^{\rm Neu}$,
     $$     
     \langle I_0 \widetilde{\mathcal{S}}I_0^*I_0 \tilde e^n, \phi_0\rangle=
     \left\langle I_0 \sum_{j=0}^\infty B_j I_0^*I_0 
     \begin{bmatrix}
\tilde e_{0,j}^n \\\tilde e_{1,j}^n
      \end{bmatrix}
     \phi_j, \phi_0\right\rangle
     =I_0B_0^{\rm Neu} I_0^*I_0 \begin{bmatrix}
\tilde e_{0,0}^n \\\tilde e_{1,0}^n
      \end{bmatrix}= -\frac12\tilde e_{0,0}^n.
     $$ 
  Furthermore, we have 
  $$
  \|I_0 \widetilde{\mathcal{S}}\tilde z\|_{L^2(\widetilde{\Omega})} \le 
  \|\widetilde{\mathcal{S}}\tilde z\|_{L^2(\widetilde{\Omega}) \times L^2(\widetilde{\Omega})}
  \leq
  C_2\|\tilde z\|_{L^2(\widetilde{\Omega}) \times L^2(\widetilde{\Omega})},
  $$
  for some constant $C_2$,
  since $\widetilde{\mathcal{S}}$ is a bounded linear map from $L^2(\widetilde{\Omega}) \times L^2(\widetilde{\Omega})$ to $L^2(\widetilde{\Omega}) \times L^2(\widetilde{\Omega})$ by Lemma~\ref{lemma::SpectralRadius}.
Then applications of the Cauchy-Schwarz, triangle inequality 
and Lemma~\ref{lemma::L2toEnergyEstimate} give
    \begin{align*}
      |\tilde e_{0,0}^{n+1}| 
      &\le 
      |\langle I_0 \widetilde{\mathcal{S}} (E P \tilde e_0^n - \tilde e_0^n), \phi_0\rangle |
           + |\langle I_0 \widetilde{\mathcal{S}}(I-I_0^*I_0) \tilde e_0^n, \phi_0\rangle|
           + |\langle I_0 \widetilde{\mathcal{S}}I_0^*I_0 \tilde e_0^n, \phi_0\rangle|
      \\
      &\le 
      C_2\|EP\tilde e^n - \tilde e^n\|_{L^2(\widetilde{\Omega}) \times L^2(\widetilde{\Omega})}
      +
      C_2\|(I-I_0^*I_0)\tilde e^n\|_{L^2(\widetilde{\Omega}) \times L^2(\widetilde{\Omega})}
      +
\frac12 |\tilde e_{0,0}^{n}|
\\&\le 
      A \|\tilde e^n\|_c + \frac12 |\tilde e_{0,0}^n|
      \le 
      A \tilde \rho^{n} + \frac12 |\tilde e_{0,0}^n|,
    \end{align*}
  where $A = C_2 \sqrt{C_1} > 0$ and $C_1 = \max\{(a-\tilde a)/\tilde c^2_\text{min},(\tilde b - b)/\tilde c^2_\text{min},1\} > 0$.

  Without loss of generality we assume $\tilde \rho \ne 1/2$ since it is possible to
  choose $\tilde a, \tilde b$ such that the problem is not at resonance with 
  $\tilde \rho \ne 1/2$. We define the sequence
    \begin{align*}
      y_n = \frac{|\tilde e_{0,0}^n|}{A} - \frac{\tilde \rho ^n}{\tilde \rho - 1/2}.
    \end{align*}
  Then
    \begin{align*}
      y_{n+1} = \frac{|\tilde e_{0,0}^{n+1}|}{A} - \frac{\tilde \rho^{n+1}}{\tilde \rho - 1/2}
      &\le 
      \frac{A \tilde \rho^n + \frac12 |\tilde e_{0,0}^n|}{A} - \frac{\tilde \rho^{n+1}}{\tilde \rho - 1/2}
      \\
      &=
      \frac12 \left(\frac{|\tilde e_{0,0}^n|}{A} - \frac{\tilde \rho^n}{\tilde \rho - 1/2}\right)
      =
      \frac12 y_n.
    \end{align*}
  Therefore $y_n \le 2^{-n} y_0$ so that
    \begin{align*}
      \lim_{n\rightarrow\infty} |\tilde e_{0,0}^n| = \lim_{n\rightarrow\infty}  A \left(y_n +  \frac{\tilde \rho^n}{\tilde \rho - 1/2}\right) = 0.
    \end{align*}
  Taking a limit of \eqref{eqn::Poincare} gives that 
  $\|\tilde e^n\|^2_{H^1(\widetilde{\Omega})\times L^2(\widetilde{\Omega})}\rightarrow 0$,
  so that we obtain convergence of the iteration in \mbox{$H^1(\Omega)\times L^2(\Omega)$}.
\end{proof}

\begin{remark}
  The above analysis is for a single spatial dimension, but we note that it in certain 
  situations it may be extended to higher dimensions. For instance, interior impedance problems 
  with constant coefficients and simple geometries
  may be extended by an appropriate enclosing box from which the above arguments 
  can give convergence. In general, it is difficult to prove convergence in higher dimensions
  in this way as  care needs to be taken to make appropriate wavespeed extensions that avoid reflections due to potentially discontinuous wavespeeds close to boundaries with impedance conditions.
\end{remark}

\section{Damped Wave/Helmholtz Equation}
As mentioned in the introduction, a popular preconditioning approach 
for solving Helmholtz problems is to introduce a damping term as in the shifted Laplacian
preconditioners \cite{erlangga2008advances}. In this
section we similarly consider the complex-valued damped wave equation
  \begin{align*}
    w_{tt} + \eta w_t = \nabla \cdot \left[c^2(x) \nabla w\right] - f(x) e^{-i \omega t},
    \qquad \eta>0,
  \end{align*}
for which we note that if $w(t,x) = u(x)e^{-i\omega t}$ then 
  \begin{align*}
    \nabla \cdot \left[c^2(x) \nabla u\right] + \left(\omega^2 + i \eta \omega\right)u = f(x),
  \end{align*}
so that we essentially have added a purely imaginary shift of 
the Laplacian 
  \begin{align*}
    \mathcal{L} = -\nabla \cdot \left[c^2(x) \nabla \right] - i\eta\omega.
  \end{align*}
Here we consider only problems with energy conserving boundary
conditions (i.e. Dirichlet or Neumann), and as a result of 
the imaginary shift of the Laplacian we note
there are no longer resonant frequencies.
While for the sake of simplicity we consider the complex-valued problem
in this section, in practice we solve the real-valued problem as presented
in Section~\ref{sec::GeneralIteration} with the filter \eqref{eq::filterstep}.
For the above complex-valued problem, we may then similarly prove an analogous 
result to Theorem~\ref{thm::FirstOrderConv}
\begin{theorem}
  The iteration \eqref{eqn::Iteration} with the complex-valued filter
  \begin{align}\label{eqn::damped_Iter}
  {\Pi} \left[
  \begin{array}{c}
  v_0\\
  v_1 
  \end{array}
  \right] = \frac{1}{T}\int_0^{T} e^{i \omega t} \left[
  \begin{array}{c}
  w(t,x)\\
  w_t(t,x) 
  \end{array}
  \right] dt,\quad
  T=\frac{2\pi}{\omega},
  \end{align}
  converges for every $\eta > 0$ with a convergence rate
  bounded by $2(1 - e^{-\eta T/2})/\eta T$.
\end{theorem}
\begin{proof}
  Suppose $(\lambda_j^2,\phi_j)$ are the eigenmodes
  of the real-valued Laplacian in the domain $\Omega$.
  We note that the shifted Laplacian now 
  has a spectrum that is $\lambda_j^2 - i \eta \omega$. 
  Expanding in terms of this basis and taking inner products,
  we can see that 
    \begin{align*}
      (\omega^2 + i \eta\omega -\lambda_j^2)u_j = f_j,
      \end{align*}
  where we expand the real and imaginary parts of $u$ and $f$ 
  as $u_j = u_j^R + i u_j^I$ and $f_j = f_j^R + i f_j^I$.
  Let the damped wave equation solution have the form
    \begin{align*}
      \sum_{n = 0}^\infty w_j(t) \phi_j(x).
    \end{align*}
  Defining $\alpha_j = \sqrt{ 4\lambda_j^2 - \eta^2}/2$,
  then the solution can be shown to be given by
    \begin{align*}
      w_j(t) = &u_j\left( e^{-i\omega t} - e^{-\frac{\eta t}{2}}\left[ \cos(\alpha_j t) +\frac{\eta - 2 i \omega}{2 \alpha_j} \sin(\alpha_j t)\right]\right) 
      \\
     &+  e^{-\frac{\eta t}{2}}\left(v_{0,j}  \cos(\alpha_j t) + \frac{\eta\sin(\alpha_j t)}{2\alpha_j}  + \frac{v_{1,j}\sin(\alpha_j t)}{\alpha_j} \right),
    \end{align*}
  from which we note that we arrive at exactly the same set of coefficients as in 
  the previous analysis if $\eta = 0$ and the real part of the solution is taken. 
  Using the complex-valued filters 
  \[
  \hat \beta(\alpha)
  := \frac{1}{T}\int_0^Te^{(i \omega - \eta/2)t}\cos(\alpha t) dt, \quad 
  \hat \gamma(\alpha)
  := \frac{1}{T}\int_0^Te^{(i \omega - \eta/2)t}\sin(\alpha t) dt,
   \]
  we can write the iteration as
  \begin{align}
    \begin{pmatrix}
      v_{0,j}^{n+1}\\v_{1,j}^{n+1}
    \end{pmatrix}
    =
    \Pi
    \begin{pmatrix}
      v_{0,j}^{n}\\v_{1,j}^{n}
    \end{pmatrix}
    =
    \left(I - \hat B_j\right)
    \begin{pmatrix}
      u_j\\i\omega u_j
    \end{pmatrix}
    +
    \hat B_j
    \begin{pmatrix}
      v_{0,j}^{n}\\v_{1,j}^{n}
    \end{pmatrix},
  \end{align}
  where if $\hat\beta_j = \hat\beta(\alpha_j)$ and $\hat\gamma_j = \hat\gamma(\alpha_j)$ then
  \begin{align*}
     \hat B_j =  
     \begin{pmatrix}
      \hat\beta_j + \frac{\eta}{2\alpha_j} \hat\gamma_j & \hat\gamma_j/\alpha_j \\
      -(\alpha_j + \frac{\eta^2}{4\alpha_j})\hat\gamma_j & \hat\beta_j - \frac{\eta}{2\alpha_j} \hat\gamma_j\\
     \end{pmatrix}.
  \end{align*}
  As in the previous analysis, we require that the spectral radius of $\hat B_j$ be less than one. The eigenvalues are given by $\hat\mu_j = \hat\beta_j \pm i \hat\gamma_j$
  so that by definition
    \begin{align}\label{eqn::damped_eig}
      |\hat\mu_j| = |\hat\beta_j \pm i \hat\gamma_j| = \left|\frac{1}{T} \int_0^T e^{i (\omega \pm \alpha_j)t}e^{-\eta t/2}\, dt \right | \le \frac{2}{\eta T} (1 - e^{-\eta T/2}) < 1,
    \end{align}
  given that $\eta > 0$.
\end{proof} 
Thus the iteration always converges in the damped case without extra
conditions on the eigenvalues.
From \eqref{eqn::damped_eig} we see that for a desired fixed rate of convergence the damping parameter $\eta$ must grow 
proportionally to $\omega$ since $\eta T\sim \eta/\omega$, and that frequency-independent convergence 
is achieved by choosing $\eta = \mathcal{O}(\omega)$.
\begin{remark}
  We note that in this section we use the complex-valued 
  filter $e^{i\omega t}/T$ instead of the usual filter, $2(\cos(\omega t) - 1/4)/T$.
  The choice of filter in \eqref{eqn::damped_Iter}, as well as performing
  the analysis using complex arithmetic, was done for the sake of 
  simplicity. The choice of filter need not be restricted to 
  $2(\cos(\omega t) - 1/4)/T$, we refer the reader to Section 2.3 and 4.1.4
  of \cite{WaveHoltz} for futher discussion on the choice of filter.
\end{remark}


\section{Analysis of Higher Order Time-Stepping Schemes for the Discrete Iteration}\label{sec:discrete}
We introduce the temporal grid points $t_n=n\Delta t$
and a spatial grid with $N$ points
together with the 
vector $w^n\in {\mathbb R}^N$ containing
the grid function values of the approximation at $t=t_n$.
We also let $f\in {\mathbb R}^N$ hold the corresponding values
of the right hand side.
The discretization of the
continuous spatial operator $-\nabla\cdot(c^2(x)\nabla)$, 
including the boundary conditions,
is denoted 
$L_h$ and it can be represented as an $N\times N$ matrix.
The values \mbox{$-\nabla\cdot(c^2(x)\nabla w)$} are then approximated by $L_hw^n$.
As in the continuous case, 
we assume $L_h$ has
the eigenmodes $(\lambda_j^2,\phi_j)$,
such that $L_h\phi_j=\lambda_j^2\phi_j$
for $j=1,\ldots,N$, where all $\lambda_j$ are real, strictly positive and
ordered as
$0\leq \lambda_1\leq\ldots\leq \lambda_N$.

We let the discrete Helmholtz solution $u$ be defined through
$$
  -L_hu + \omega^2 u = f.
$$
The numerical approximation of the iteration operator is
denoted $\Pi_h$, and it
is implemented as follows. Given $v\in\mathbb{R}^N$, we use
the leap frog method to solve the wave equation and add in 
higher order corrections as in the Modified Equation (ME) approach 
\cite{shubin1987modified,anne2000construction}.
For a general $2m$ scheme, recall that via Taylor expansion
  \begin{align*}
    \frac{w^{n+1} - 2w^n+w^{n-1}}{\Delta t^2} = w_{tt} + 2 \sum_{k=2}^\infty \frac{\Delta t^{2(k-1)}}{(2k)!} \frac{\partial^{2k}}{\partial t^{2k}} w^n.
  \end{align*}
Then using the PDE to convert time derivatives to spatial derivatives we get the expression
  \begin{align*}
    \frac{\partial^{2k} }{\partial t^{2k}} w^n \approx L_h^k w^n + \cos(\omega t_n) \sum_{\ell=0}^{k-1} (-1)^{k+\ell} \omega^{2(k - \ell - 1)}L_h^\ell f,
  \end{align*}
for $k = 1,2, \dots$. Then for a $2m$ order scheme we have
  \begin{align}\label{eqn::MEstepper}
    \frac{w^{n+1} - 2w^n+w^{n-1}}{\Delta t^2} - 2 \sum_{k=2}^m \frac{\Delta t^{2k-2}}{(2k)!} \left[L_h^k w^n + \cos(\omega t_n) \sum_{\ell=0}^{k-1} (-1)^{k+\ell} \omega^{2(k - \ell - 1)}L_h^\ell f\right] \\
    = L_h w^n - f\cos(\omega t_n),\nonumber
  \end{align}
with time-step $\Delta t=T/M$ for some integer $M$, and initial data 
$$
   w^0 = v,\qquad w^{-1} = v + \sum_{k=1}^{m} \frac{(-1)^k\Delta t^{2k}}{(2k)!} \left[-L_h^k v + \sum_{\ell=0}^{k-1} (-1)^{\ell} \omega^{2(k - \ell - 1)}L_h^\ell f\right].
$$
The trapezoidal rule is then used to compute $\Pi_h v$,
\be\lbeq{trapzrule}
\Pi_h v
= 
   \frac{2\Delta t}{T}
   \sum_{n=0}^M \eta_n \left(\cos(\omega t_n)-\frac14\right)w^n,
   \qquad \eta_n = \begin{cases}
   \frac12,& \text{$n=0$ or $n=M$},\\
   1, & 0<n<M.
   \end{cases}
\ee
We may then prove the following theorem
that is a generalization of Theorem 2.4 
of \cite{WaveHoltz}.
\begin{theorem}\label{thm::DiscreteTS}
Suppose that $L_h$ has
real and strictly positive eigenvalues $\lambda_j$ and that
there are no resonances, such that
$\delta_h=\min_j|\lambda_j-\omega|/\omega>0$.
Moreover, assume that $\Delta t$
satisfies the stability and accuracy requirements
\be\lbeq{cfl}
  \Delta t< \frac{2}{\lambda_N+2\omega/\pi},
  \qquad \Delta t\omega\leq \min(\delta_h,1).
\ee
Then the fixed point iteration $v^{(k+1)}=\Pi_h v^{(k)}$ with
$v^{(0)}=0$ converges to $v^\infty$ which 
is a solution to the discretized Helmholtz equation,
  \begin{align*}
      -L_hv^\infty + \tilde\omega^2 v^\infty = f,
  \end{align*}
with the modified frequency $\tilde\omega$,
defined as the smallest positive real number satisfying
  \begin{align*}
      \sin^2(\omega\Delta t/2) = \sum_{j=1}^m \frac{(-1)^{j+1}\left(\Delta t \tilde \omega\right)^{2j}}{2(2j)!}, 
  \end{align*}
where $2m$ is the order of the ME time-stepping scheme.
Moreover, there are constants $C_m<1$ and $C_m'$ only depending om $m$ such that
$$
|\omega - \tilde\omega| \leq C_m\Delta t^{2m}\omega^{2m+1},
\qquad \|u - v^\infty\|_2 \leq 
C_m'\Delta t^{2m}\omega^{2m}\delta_h^{-2}\|f\|_2.
$$
The convergence rate is at least $\rho_h=\max(1-0.3\delta_h^2,0.6)$.
\end{theorem} 

\begin{proof}
We expand all functions in eigenmodes of $L_h$,
\begin{gather*}
  w^n = \sum_{j=1}^N w_j^n \phi_j,
  \qquad
  f = \sum_{j=1}^N f_j \phi_j,\qquad
  u = \sum_{j=1}^N u_j \phi_j,\qquad
  \\
  v = \sum_{j=1}^N v_j \phi_j,\qquad
  v^\infty = \sum_{j=1}^N v^\infty_j \phi_j.
\end{gather*}
Then the Helmholtz eigenmodes of $u$ and $v^\infty$ satisfy
$$
u_j = \frac{f_j}{\omega^2-\lambda_j^2},\qquad
v^\infty_j = \frac{f_j}{\tilde\omega^2-\lambda_j^2}.
$$
We note that $\tilde\omega$ is well-defined
by Lemma \ref{lem::welldefined} in Appendix~\ref{sec:well-defined}.
The same lemma also shows the bound on $|\omega-\tilde{\omega}|$,
which implies that $\tilde\omega$ is not resonant and
$v^\infty_j$ is well-defined
for all $j$, since by
\eqref{eqn::Freq_error} and \eq{cfl}
\begin{align*}
|\tilde\omega-\lambda_j|\geq |\omega-\lambda_j|-
|\tilde\omega-\omega|\geq \omega\delta_h - 
C_m\Delta t^{2m}\omega^{2m+1}
&\geq
\omega\left(\delta_h - C_m\min(\delta_h,1)^{2m}\right)
\\
&\geq
\omega\delta_h(1 - C_m)>0.  
\end{align*}
The wave solution eigenmodes to \eqref{eqn::MEstepper} are given by the difference equation
\begin{align}\label{eqn::DifferenceEQ}
  w_j^{n+1} - 2w_j^n+w_j^{n-1} + &2 \left[\sum_{k=1}^m \frac{(-1)^{k+1}\Delta t^{2k}\lambda_j^{2k}}{(2k)!} \right]w_j^n 
  \\
  &= 2\left[\sum_{k=1}^m  \frac{(-1)^{k}\Delta t^{2k}}{(2k)!}\sum_{\ell=0}^{k-1} \omega^{2(k - \ell - 1)}\lambda_j^{2\ell} \right]f_j\cos(\omega t_n),\nonumber
\end{align}
with initial data
\begin{gather*}
w_j^{0} =
v_j,\\
w_j^{-1} =
v_j\left(1 + \sum_{k=1}^m\frac{(-1)^k\Delta t^{2k}}{(2k)!} \lambda_j^{2k}\right) + f_j\left(\sum_{k=1}^{m}\frac{(-1)^k\Delta t^{2k}}{(2k)!}\sum_{\ell=0}^{k-1} \omega^{2(k - \ell - 1)}\lambda_j^{2\ell}\right).
\end{gather*}
By \eq{cfl}, the discrete solution is stable and given by
\begin{align}\label{eqn::discreteSoln}
     w^n_j = 
   (v_j-v^\infty_j)\cos(\tilde\lambda_j t_n) + v^\infty_j\cos(\omega t_n),
\end{align}
where $\tilde \lambda_j$ is well-defined, by 
\eqref{eqn::lambound},
as the smallest positive real number satisfying
  \begin{align}\label{eq:freq_shift}
    \sin^2(\tilde\lambda_j\Delta t/2) = \sum_{k=1}^m \frac{(-1)^{k+1}\left(\Delta t \lambda_j\right)^{2k}}{2(2k)!}.
  \end{align}

For $m\ge 2$, we have that 
$|\omega - \tilde \omega| \le C_m\Delta t^{2m}\omega^{2m+1} \le \Delta t^2 \omega^3/24$
since $C_m = 5/(2m+2)! \le 1/24$ by Lemma~\ref{lem::welldefined}. We may then apply the 
following lemma, restated from \cite{WaveHoltz}, to obtain convergence of the discrete iteration 
(we note that the proof of Lemma~\ref{filterlemma2} requires a simple modification for the case $m=1$ 
and is thus not presented here for the sake of brevity).

\begin{lemma}\label{filterlemma2}
Under the assumptions of Theorem~\ref{thm::DiscreteTS},
\begin{equation}\lbeq{betahest}
   \max_{1\leq j\leq N} |\beta_h(\tilde\lambda_j)|\leq \rho_h =:\max(1-0.3\delta_h^2,0.6).
\end{equation}
\end{lemma}

Letting $e = u - v^\infty$ be the error in the discrete solutions, 
the components of the error in the basis of the Laplacian satisfy
  \begin{align*}
    |e_j| = |u_j - v_j^\infty| 
    &= \left|f_j\left(\frac{1}{\omega^2 - \lambda_j^2} - \frac{1}{\tilde\omega^2 - \lambda_j^2}\right)\right| 
    \\
    &= \left|f_j(\tilde\omega - \omega)\right|
    \left|\frac{\omega+\tilde\omega}
    {(\omega+\lambda_j)(\tilde\omega+\lambda_j)
    (\omega-\lambda_j)(\tilde\omega-\lambda_j)}
    \right|
    \\
    &\leq C_m|f_j|\Delta t^{2m}\omega^{2m+1}
    \left|\frac{\omega^{-1}+\tilde\omega^{-1}}
    {(1-C_m)\delta_h^2}
    \right|
    \\
    &\leq \frac{C_m(2-C_m)}{(1-C_m)^2}
    |f_j|\Delta t^{2m}\omega^{2m}\delta_h^{-2}
    =:C_m'
    |f_j|\Delta t^{2m}\omega^{2m}\delta_h^{-2},
  \end{align*}
  where we also used the fact that
  $$
    \frac{\omega}{\tilde{\omega}}
    =\frac{\omega}{\omega+\tilde{\omega}-\omega}
    \leq
    \frac{\omega}{\omega-C_m\Delta t^{2m}\omega^{2m+1}}
   = \frac{1}{1-C_m\Delta t^{2m}\omega^{2m}}
   \leq \frac{1}{1-C_m}.
  $$
This gives
  \begin{align*}
     \|u - v^\infty\|_2 = \|e\|_2 
     \leq C_m'\Delta t^{2m}\omega^{2m}\delta_h^{-2}\|f\|_2.
  \end{align*}
concluding the proof of the theorem.  
\end{proof}

\begin{remark}
  As alluded to in Remark 6 of \cite{WaveHoltz}, knowledge of how a particular 
  discretization approximates the eigenvalues of the continuous operator can be 
  used to improve the iteration. In fact, the above error due to time 
  discretization can be removed by defining $\bar \omega$ by the relation
    \begin{align*}
      \sin^2(\bar\omega\Delta t/2) = \sum_{k=1}^m \frac{(-1)^{k+1}\left(\Delta t \omega\right)^{2k}}{2(2k)!}.
    \end{align*}
  Then using $f \cos(\bar \omega t_n)$ instead of $f \cos(\omega t_n)$ in the time-stepping \eqref{eqn::MEstepper}, 
  in addition to the modified trapezoidal quadrature rule (first introduced in \cite{peng2021emwaveholtz})
    \begin{align}
    \Pi_h v
    = 
       \frac{2\Delta t}{T}
       \sum_{n=0}^M \eta_n\frac{\cos(\omega t_n)}{\cos(\bar \omega t_n)} \left(\cos(\omega t_n)-\frac14\right)w^n,
       \qquad \eta_n = \begin{cases}
       \frac12,& \text{$n=0$ or $n=M$},\\
       1, & 0<n<M,
       \end{cases}\label{eqn::quadrature}
    \end{align}
  gives that the limit will be precisely the discrete Helmholtz solution, $v^\infty = u$, as 
  long as the time-step size is chosen so that $\cos(\bar \omega t_n) \ne 0$.
  Moreover, the first time-step restriction of \eq{cfl} arising from the usual 
  ${\rm CFL}$ condition for the second order scheme may be relaxed (expressions for which 
  may be found in \cite{Gilbert2008}) though the condition 
  $\Delta t\omega \leq \min(\delta_h,1)$ may be more restrictive for 
  problems close to resonance. We additionally note that 
  in \cite{stolk2020timedomain} an alternative approach to remove time-discretization
  error was presented, however the approach modified the time-stepping scheme whereas we
  modify the frequency of the forcing and update our quadrature rule.
\end{remark}



\section{Wave Equation Solvers}
In this section we briefly outline the numerical methods we use in the experimental section below. We consider both discontinuous Galerkin finite element solvers and finite difference solvers. In all the experiments we always use the trapezoidal rule to compute the integral in the WaveHoltz iteration. 

\subsection{The Energy Based Discontinuous Galerkin Method} \label{sec:dG}
Our spatial discretization is a direct application of the formulation described for general second order wave equations in \cite{Upwind2,el_dg_dath}. Here we outline the spatial discretization for the special case of the scalar wave equation in one dimension and refer the reader to \cite{Upwind2} for the general case. 

The energy of the scalar wave equation is 
\[
H(t) = \int_{D} \frac{v^2}{2} + G(x,w_x) dx,
\]
where 
\[
 G(x,w_x) = \frac{c^2(x)w_x^2}{2},
\] 
is the potential energy density, $v$ is the velocity (not to be confused with the iterates $v^n$ above) or the time derivative of the displacement, $v = w_t$. The wave equation, written as a second order equation in space and first order in time then takes the form 
\begin{eqnarray*}
w_t &=& v, \\
v_t &=& - \delta G, 
\end{eqnarray*}
where $\delta G$ is the variational derivative of the potential energy 
\[
\delta G = - (G_{w_x})_x = -(c^2(x)w_x)_x.
\]
For the continuous problem the change in energy is   
\begin{equation}
\f{d H(t)}{dt} = \int_{D} v v_t + w_t  (c^2(x)w_x)_x \,dx = [ w_t  (c^2(x)w_x)]_{\partial D}, \label{eq:energy_derivative}
\end{equation}
where the last equality follows from integration by parts together with the wave equation. Now, a variational formulation that mimics the above energy identity can be obtained if the equation $v-w_t=0$ is tested with the variational derivative of the potential energy. Let $\Omega_j$ be an element and $\Pi^s(\Omega_j)$ be the space of polynomials of degree $s$, then the variational formulation on that element is:
\begin{prob}
Find $v^h \in \Pi^s(\Omega_j)$, $w^h \in \Pi^{r}(\Omega_j)$ such that for all
$\psi \in \Pi^s(\Omega_j)$, $\phi \in \Pi^{r}(\Omega_j)$ 
\begin{eqnarray}
\int_{\Omega_j} c^2 \phi_x \left( \f {\pa w^h_x}{\pa t}-v^h_x \right) dx & = & 
[c^2\phi_x \cdot n \left( v^{\ast}-v^h \right)]_{\pa \Omega_j}, 
\label{var1} \\
\int_{\Omega_j} \psi \f {\pa v^h}{\pa t} + c^2 \psi_x \cdot  w^h_x \, dx& = &
 [\psi \, (c^2\,w_x)^{\ast}]_{\pa \Omega_j}. \label{var2} 
\end{eqnarray}
\end{prob}

Let $[[f]]$ and $\{f\}$ denote the jump and average of a quantity $f$ at the interface between two elements, then, choosing the numerical fluxes as 
\begin{eqnarray*}
v^{\ast}  &=& \{v\} -\tau_1 [[c^2\,w_x]]\\
(c^2\,w_x)^{\ast} &=& \{ c^2\,w_x \}  -\tau_2 [[v]],
\end{eqnarray*}
will yields a contribution $ -\tau_1 ([[c^2\,w_x]])^2 -\tau_2 ([[v]])^2$ from each element face to the change of the discrete energy
\[
\f{d H^h(t)}{dt} = \frac{d}{dt} \sum_{j} \int_{\Omega_j} \frac{(v^h)^2}{2} + G(x,w^h_x).
\]
Physical boundary conditions can also be handled by appropriate specification of the numerical fluxes, see \cite{Upwind2} for details. The above variational formulation and choice of numerical fluxes results in an energy identity similar to (\ref{eq:energy_derivative}). However, as the energy is invariant to certain transformations the variational problem does not fully determine the time derivatives of $w^h$ on each element and independent equations must be introduced. In this case there is one invariant and an independent equation is $\int_{\Omega_j} \left( \f {\pa w^h}{\pa t}-v^h\right) = 0$. For the general case and for the elastic wave equation see \cite{Upwind2} and  \cite{el_dg_dath}.  

In this paper we always choose $\tau_i > 0$ (so-called upwind or Sommerfeld fluxes) and we always choose the approximation spaces to be of the same degree $r=s$. These choices result in methods that are $r+1$ order accurate in space. 

\subsection{Symmetric Interior Penalty Discontinuous Galerkin Method}
In addition to the above energy DG method, we also consider the Symmetric Interior Penalty DG (SIPDG) discretization, 
\cite{GSSwave}, for examples in two dimensions. The bilinear form in this case is 
	\begin{align*}
		a_h(u,v) =  \sum_{K \in \mathcal{T}_h} \int_K c^2 \nabla u \cdot \nabla v \, dx 
				- \sum_{f \in \mathcal{F}_h} & \int_F [[u]]\cdot \{c^2 \nabla v\}  
				- [[v]]\cdot \{c^2 \nabla u\} 
				\\
				&+ \gamma h_F^{-1} c^2[[u]]\cdot [[v]] \, ds,
	\end{align*}
where $\mathcal{T}_h$ is a collection of triangular elements, 
$\mathcal{F}_h$ is the collection of element faces, 
$h_F$ is the diameter of the edge or face $F$, and $\gamma$ is the interior
penalty stabilization parameter which must be chosen to be sufficiently
large to ensure the system is positive-definite.

\subsection{Finite Difference Discretizations} \label{sec:FD}
For the finite difference examples in a single dimension, we consider discretizations 
by uniform grids $x_i = x_L+ih_x,$, with $i = -1,\ldots,n+1$ and $h_x = (x_R-x_L) / n$.
To impose impedance boundary conditions of the form $w_t \pm \vec{n} \cdot \nabla w  = 0$ 
we evolve the wave equation as a first order system in time according to the semi-discrete approximation 
\begin{align*}
\frac{d v_{i}(t)}{dt}  =  (D_+D_-) w_{i}, \qquad
\frac{d w_{i}(t)}{dt}  = v_{i},   
\end{align*}
and for the boundaries we find the ghost point values by enforcing 
\begin{equation}
v_{0} - D_0 w_{0} = 0, \quad v_{n} - D_0 w_{n} = 0.
\end{equation}
Here we have used the standard forward, backward and centered finite difference operators, 
for example $hD_+w_i = w_{i+1}-w{i}$ etc.

\subsection{Time Discretization}
For some of the numerical examples in a single dimension, we use either an explicit second order accurate centered discretization of $w_{tt}$ or use the higher order corrected ME methods described in Section~\ref{sec:discrete}.      

For the DG discretizations we employ Taylor series time-stepping in order to match the order of accuracy in space and time.  Assuming that all the degrees of freedom have been assembled into a vector ${\bf w}$ we can write the semi-discrete method as ${\bf w}_t = Q {\bf w} $ with $Q$ being a matrix  representing the spatial discretization. Assuming we know the discrete solution at the time $t_n$ we can advance it to the next time-step $t_{n+1} = t_n + \Delta t$ by the simple formula
\begin{align*}
{\bf w}(t_{n}+\Delta t) &= {\bf w}(t_{n}) + \Delta t {\bf w}_t(t_{n}) +  \frac{(\Delta t)^2}{2!}{\bf w}_{tt}(t_{n}) \ldots
 \\&= {\bf w}(t_{n}) + \Delta t Q {\bf w}(t_{n}) +  \frac{(\Delta t)^2}{2!} Q^2 {\bf w}(t_{n}) \ldots.
\end{align*} 
The stability domain of the Taylor series which truncates at time derivative number $N_{\rm T}$ includes the imaginary axis if $ {\rm mod} (N_{\rm T},4) = 3$ or ${\rm mod} (N_{\rm T},4) = 0$. However as we use a slightly dissipative spatial discretization the spectrum of our discrete operator will be contained in the stability domain of all sufficiently large choices of $N_{\rm T}$ (i.e. the $N_{\rm T}$ should not be smaller than the spatial order of approximation).

\section{Numerical Examples}\label{sec:numexp}
In this section we illustrate the properties of the proposed iteration and its Krylov accelerated 
version by a sequence of numerical experiments in one and two spatial dimensions.

\subsection{Examples in One Dimension}

\subsubsection{Convergence Rate for Impedance Boundary Conditions}
In \cite{WaveHoltz}, an application of Weyl asymptotics \cite{weyl:11} revealed that the mininal relative gap to resonance, $ \delta = \min_j |\omega - \lambda_j|/\omega$ where $\lambda_j^2$ are the eigenvalues of the Laplacian, shrinks as $\omega^{-d}$ where $d$ is the spatial dimension of the Helmholtz problem of interest. Analysis of the symmetric, positive definite formulation of the iteration then yielded a convergence rate of $1 - \mathcal{O}(\delta^2) \approx 1 - \mathcal{O}(\omega^{-2d})$. However, numerical experiments with Helmholtz problems with certain open/outflow boundary conditions suggest a much more attractive convergence rate than the unacceptable $1 - \mathcal{O}(\omega^{-2d})$ rate. A natural question then is whether or not this seemingly pessimistic convergence rate can be observed for outflow boundary conditions which are much more common in practical applications.

\begin{figure}[htpb]
\graphicspath{{figures/}}
\begin{center}
\includegraphics[width=0.31\textwidth]{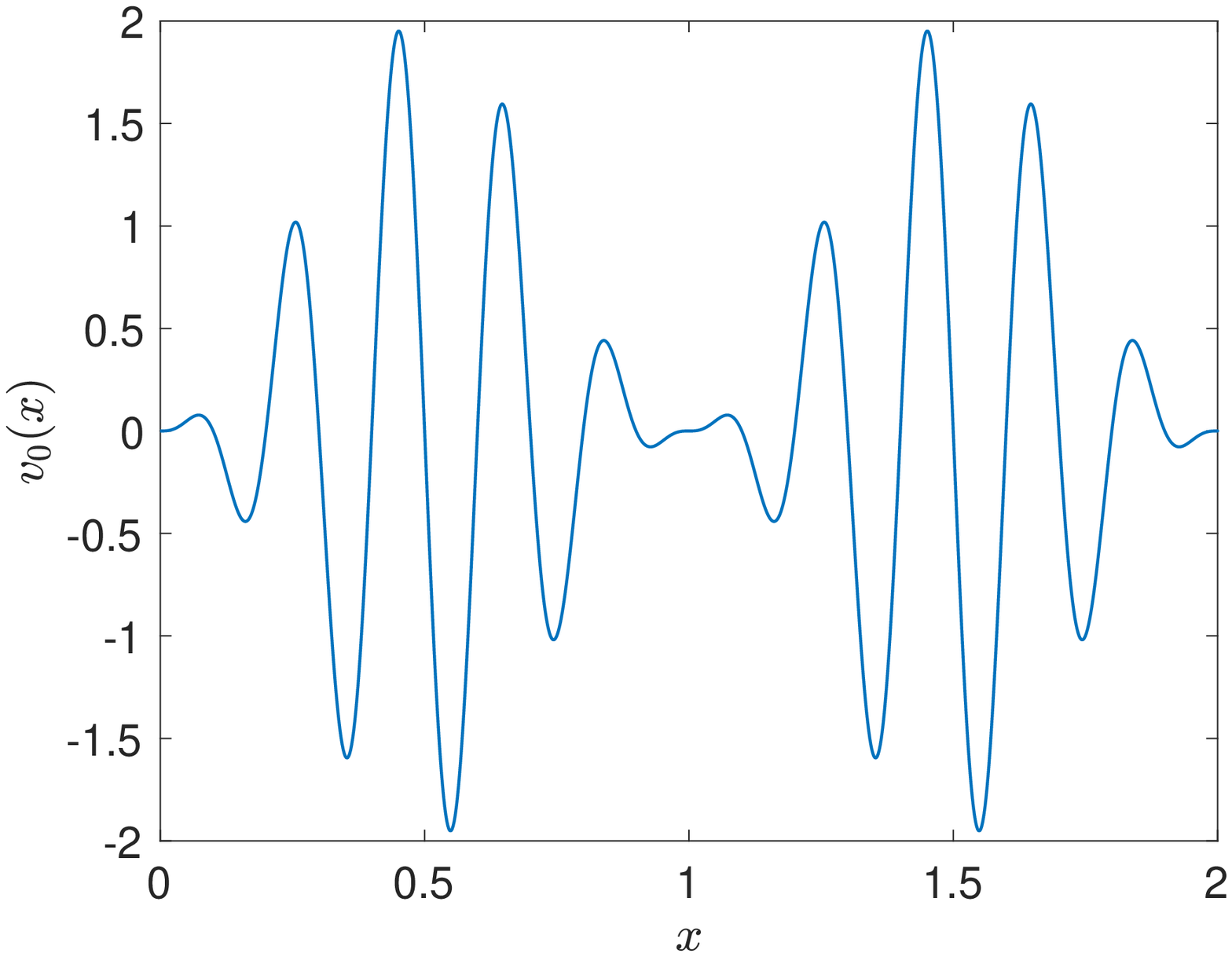}
\includegraphics[width=0.32\textwidth]{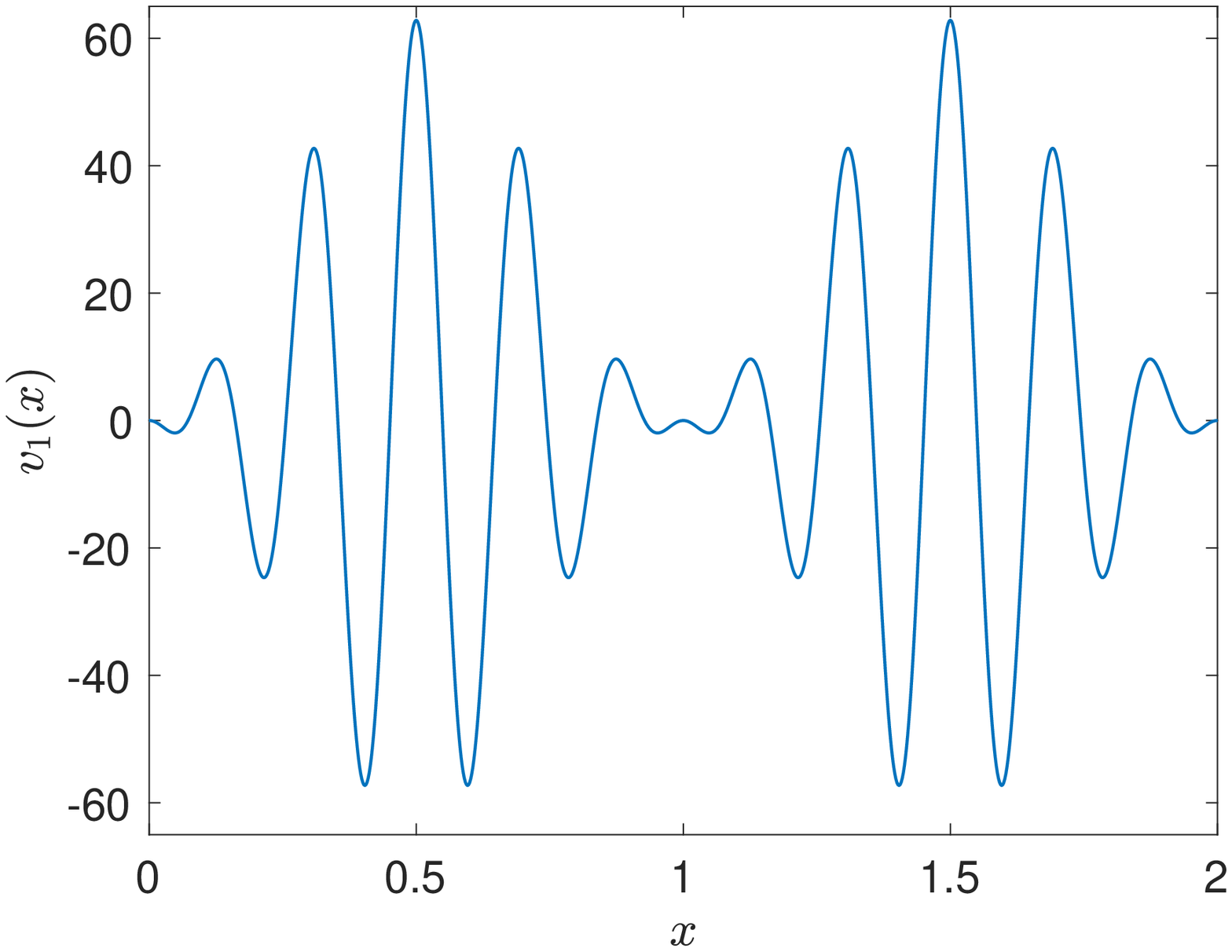}
\includegraphics[width=0.32\textwidth]{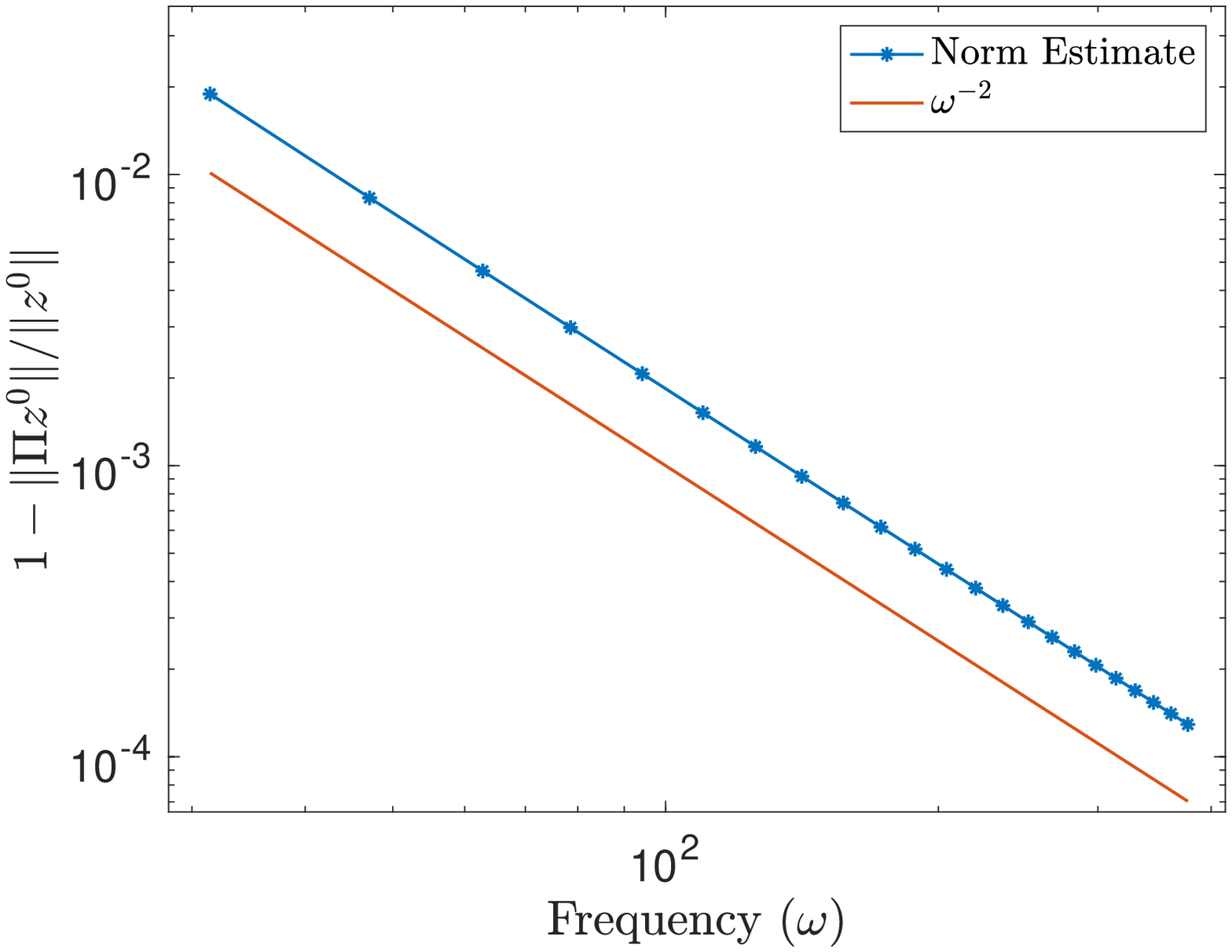}
\caption{(Left, Middle) The initial conditions $v_0$ and $v_1$ for a Helmholtz frequency of $\omega = 10 \pi$. (Right) The estimate of the quantity $1 - \|\mathcal{S}\|$ with increasing Helmholtz frequency $\omega$. \label{fig:1d_ibc_conv}}
\end{center}
\end{figure}

To that end, we consider a set of sample Helmholtz problems in a single spatial dimension with a constant (normalized) speed of sound, $c = 1$, in the domain $0 \le x \le 2$ where we impose the impedance boundary condition $w_t + \vec{n}\cdot w_x = 0$,
which we note is equivalent to the Sommerfeld radiation condition. The Helmholtz problem under consideration has no forcing and so $f = 0$. We formulate the wave equation in first order form and apply the extended iteration \eqref{eqn::Iteration} since the boundary conditions do not conserve energy. The Laplacian is discretized with a standard three-point finite difference approximation, and a fourth order Taylor scheme is used for time-stepping. 
We define the initial conditions as 
  \begin{align*}
    v_0(x) = \sin(\omega x) - \frac{1}{2}\left(\sin((\omega+2\pi) x)+\sin((\omega-2\pi) x) \right), 
    \quad 
    v_1(x) = -\frac{d}{dx} v_0(x),
  \end{align*}
which are shown in Figure~\ref{fig:1d_ibc_conv}.

By definition, $\|\mathcal{S}\|_2 = \sup_{\|z\|_2 \ne 0} \|\mathcal{S} z\|_2/\|z\|_2 \ge \|\mathcal{S} z^0\|_2/\|z^0\|_2$ so that if \newline $\|S z^0\|_2/\|z^0\| \approx 1-\mathcal{O}(\omega^{-2})$ is observed then the estimate of the spectral radius of the fixed point operator $\mathcal{S}$ is tight even for the problem with impedance boundary conditions. We consider a sweep of Helmholtz frequencies \mbox{$\omega = 10 \pi, 15\pi, 20\pi, \dots, 120 \pi$} with fifty points per wavelength and a ${\rm CFL}$ number of $10^{-1}$ for the solution of the wave equation. The results of this experiment are shown in Figure \ref{fig:1d_ibc_conv}.

On the left of Figure \ref{fig:1d_ibc_conv} we see the first part of the initial condition $v^0$ for 
a frequency of $\omega = 10\pi$. We note that this specific initial condition is constructed such that 
it is close to a resonant mode -- which the filter-transfer function $\beta$ weakly damps -- as well as 
being close to zero at the boundary so that a negligible amount of energy exits the system due to the 
impedance boundary conditions in a single iteration. These two defining characteristics of the initial
 condition lead to the norm estimate of the fixed-point iteration operator $\mathcal{S}$ on the right 
 of Figure \ref{fig:1d_ibc_conv}. We observe that the norm of $\mathcal{S}$ does indeed approach unity 
 at a rate of $\omega^{-2}$, as predicted by theory. Thus, while the preceeding analysis ``artificially'' 
 leveraged energy conserving boundary conditions to obtain an estimate of the convergence rate for open 
 problems, it is possible to realize the `worst-case' rate implied by the energy conserving regime.

\begin{figure}[H]
\graphicspath{{figures/}}
\begin{center}
\includegraphics[width=0.5\textwidth]{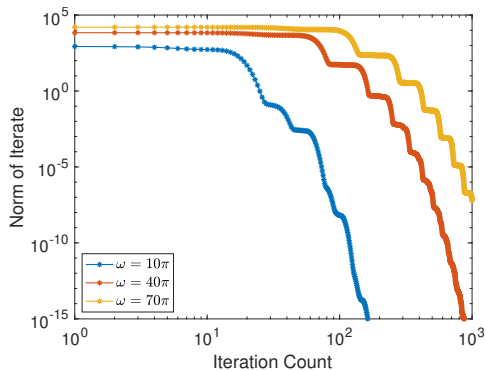}
\caption{The norm of WaveHoltz iterates for increasing Helmholtz frequencies of  $\omega = 10\pi, 40\pi,$ and $70\pi$ for the
adversarial example of Figure~\ref{fig:1d_ibc_conv}. \label{fig:PointwiseConvergence}}
\end{center}
\end{figure}

\begin{remark}
We note that the estimate for the convergence rate is a \textbf{pointwise} estimate. Repeated application
of the fixed-point iteration will (eventually) remove the modes close to resonance and a faster convergence
rate is observed. In Figure~\ref{fig:PointwiseConvergence} we repeat the above experiment for 
the frequencies $\omega = 10\pi, 40\pi,$ and $70\pi$ but continue the iteration until the iterates converge
to the zero solution. We observe that after an initial phase the rate of convergence of the iterates
to the solution increases significantly since the data has propagated and exited the domain. 
We believe that the average behavior over many fixed-point iterations 
leads to the much more attractive rates seen in the Krylov-accelerated numerical experiments of \cite{WaveHoltz}.
Moreover, this example was pathologically constructed and we note that so far we have been unable to construct
initial conditions to realize the worst-case rate in higher than one dimension.
\end{remark}

Assuming radially symmetric solutions to the Helmholtz equation, it is possible to cast higher dimensional problems as 1D problems. 
We now consider solving an analagous problem in cylindrical (2D) and spherical (3D) coordinates with
radial coordinate $r$.
We use a second order finite difference 
discretization (see \cite{morton_mayers_2005} for details) on the unit ball, $r\in[0,1]$,
with an impedance boundary condition at $r = 1$.  The initial condition is similar to the previous example,
  \begin{gather*}
    \quad v_0(r) = \sin(\omega (r+1)) - \frac{1}{2}\left(\sin((\omega+2\pi) (r+1))+\sin((\omega-2\pi) (r+1)) \right),\\
    \quad v_1(r) = -\frac{d}{dr} v_0(r).
  \end{gather*}
We consider a set of frequencies $10\pi, 11\pi, \dots, 30 \pi$ and use fifty points per wavelength in the computation with a ${\rm CFL}$ of $10^{-2}$. Below we show the results of the experiment. 
\begin{figure}[htpb]
\graphicspath{{figures/}}
\begin{center}
\includegraphics[width=0.47\textwidth]{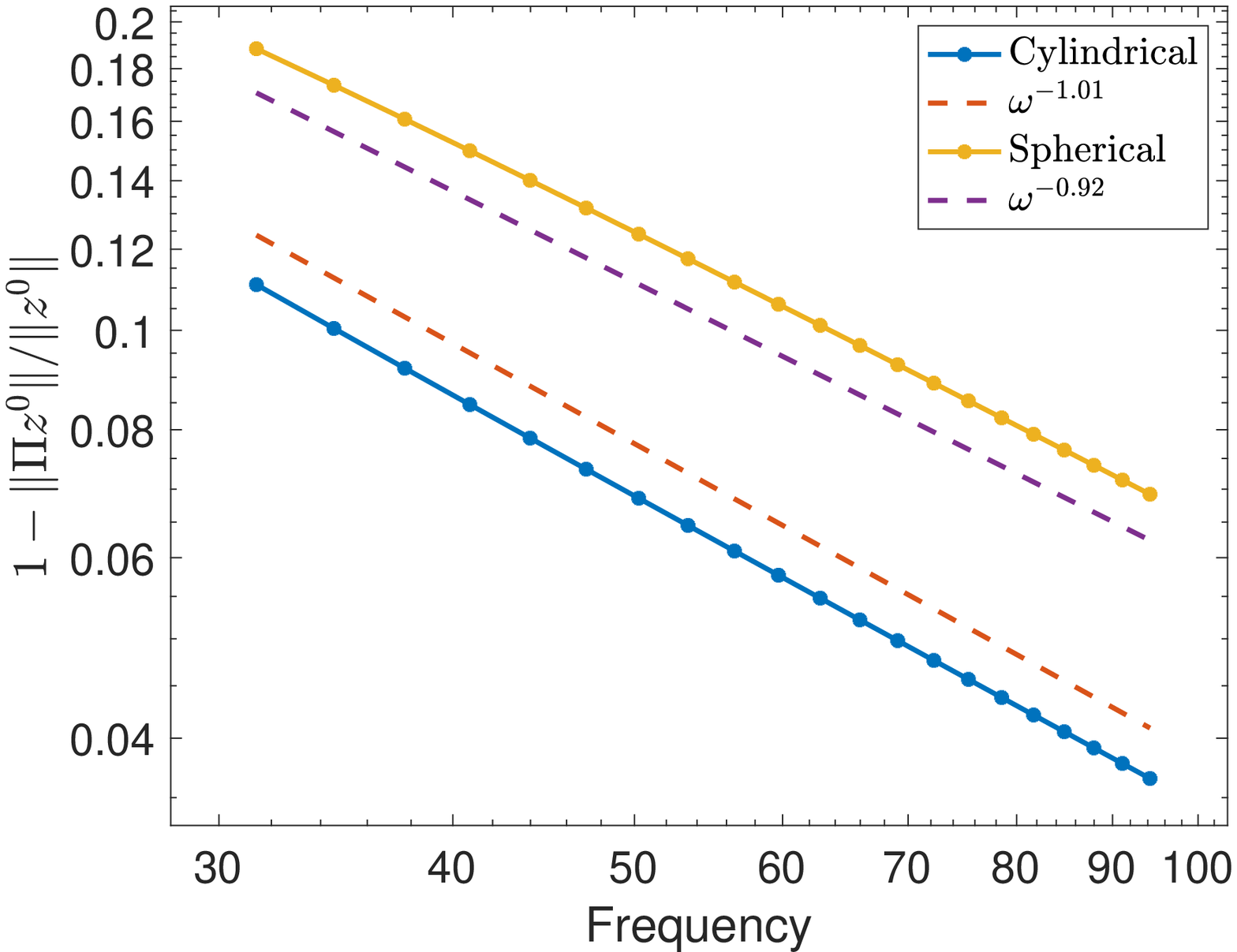}
\includegraphics[width=0.49\textwidth]{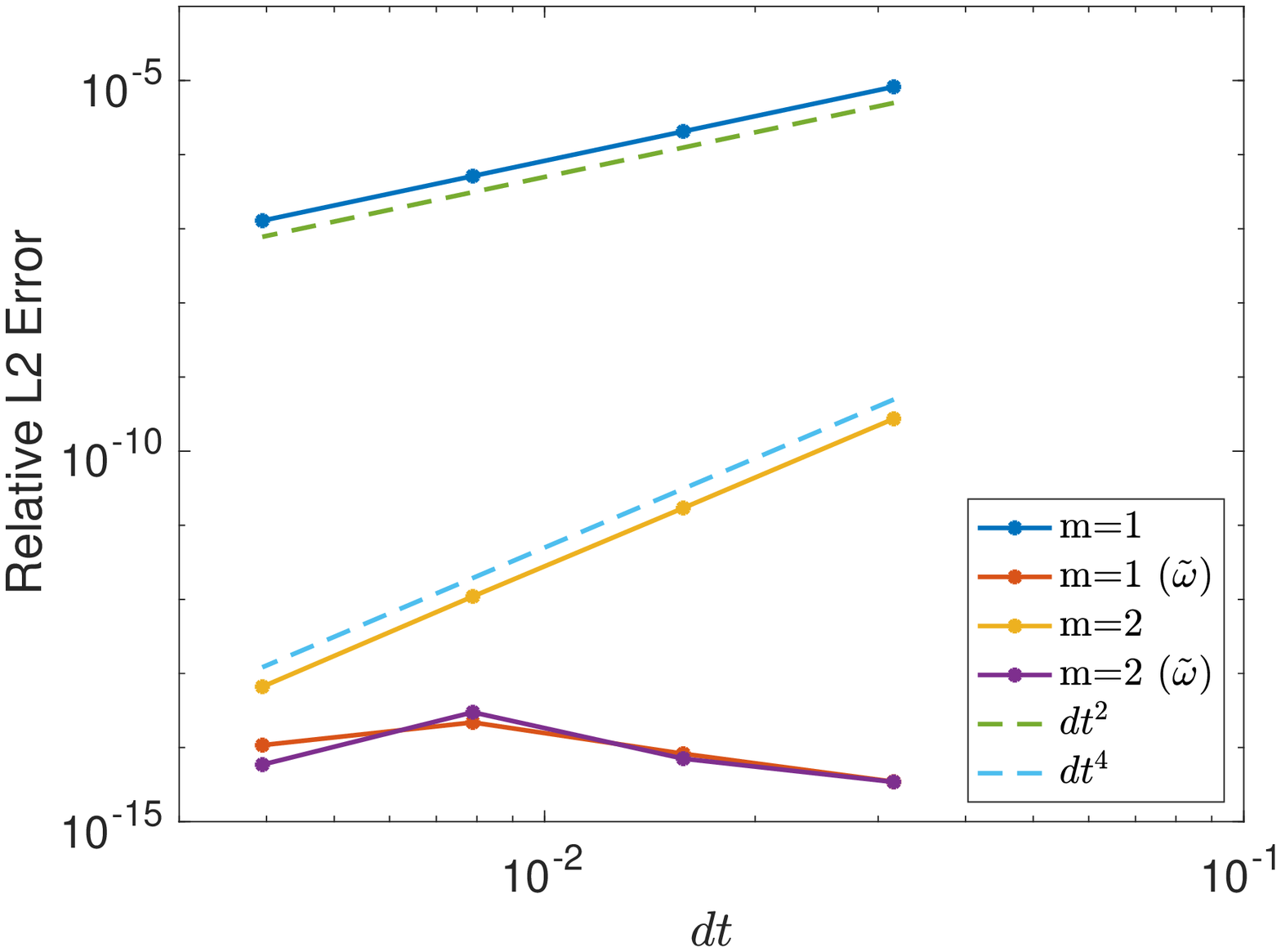}
\caption{Left: The estimate of the quantity $1 - \|\mathcal{S}\|$ with increasing Helmholtz frequency $\omega$ for a radially symmetric initial condition. Right: Convergence of the discrete WaveHoltz solution to the true solution of the discrete Helmholtz problem with fixed spatial discretization. The curves labeled with $\tilde \omega$ indicate the solution using the modified quadrature \eqref{eqn::quadrature}. \label{fig:CylSphConv}}
\end{center}
\end{figure}
From the left of Figure~\ref{fig:CylSphConv} we observe that the norm of $\mathcal{S}$ approaches unity at a nearly linear rate in the frequency $\omega$ in 2D and 
a sublinear rate for the 3D problem, both of which are more favorable than the quadratic rate in a single spatial dimension.
\begin{remark}
  From the left of Figure~\ref{fig:CylSphConv} it is clear that with a fixed discretization and initial condition, the convergence rate improves with increasing dimension.
  This is perhaps unsurprising given an increase in the local energy decay rate for the wave equation from two to three dimensions, along with 
  a richer set of directions in which waves may propagate and leave the domain.
\end{remark}

\subsubsection{Time Discretization}
We consider solving the Helmholtz equation with $c = 1$ and constant exact solution
  \begin{align*}
    u(x) = 1, \quad 0\le x \le 1.
  \end{align*}
We take the frequency to be $\omega = 1$ and consider Dirichlet boundary conditions. We discretize the Laplacian with the 
standard three-point finite difference stencil and note that there is no error (aside from truncation errors) in the solution by a direct solution of the discrete Helmholtz equation. We use a centered modified equation time-stepping scheme of both second and fourth 
order, with both the original frequency and a modified frequency $\tilde \omega$ with corresponding quadrature to remove time discretization errors. We use the WaveHoltz iteration as a fixed-point iteration with a convergence criterion that the relative $L_2$ norm between successive iterations is smaller than $10^{-13}$.
Using the original frequency in the calculation, we see on the right of Figure~\ref{fig:CylSphConv} that the WaveHoltz solution converges to the discrete Helmholtz solution with the same order as that of the time-step scheme used. With the modified frequency and quadrature, however, we see that the WaveHoltz iteration converges to the discrete Helmholtz solution up to roundoff errors.

\begin{remark}
While only centered time-stepping schemes are presented here, this approach can be 
extended to arbitrary time-steppers. A careful discrete analysis of the iteration 
isolated to a single eigenmode of the wave solution reveals what the modified frequency
should be, and a modified quadrature as outlined above removes the time discretization error 
from the converged WaveHoltz solution. Thus, the choice of a time-stepper need not need be 
restricted to have the same order as the spatial discretization. With a corrected scheme 
it may be more advantageous to take as large a time-step as possible with a low order 
time-stepper.
\end{remark}

\subsubsection{Convergence Rate for Damped Helmholtz Equations}
To study how the number of iterations scale with the Helmholtz frequency $\omega$ we solve the wave equation on the domain $x\in[-6,6]$ with constant wave speed $c^2(x) = 1$ and with a forcing 
\begin{align*}
  f(x) = \omega^2 e^{-(\omega x)^2},
\end{align*}
that results in the solution being $\mathcal{O}(1)$ for all $\omega$. We discretize using the energy based DG method discussed above and use upwind fluxes which adds a small amount of dissipation. We keep the number of degrees of freedom per wave length fixed by letting the number of elements be $5 \lceil \omega \rceil$. We always take the polynomial degree to be 7, the number of Taylor series terms in the time-stepping to be 6, and use WHI accelerated by GMRES without restarts.

\begin{figure}[htpb]
\graphicspath{{figures/}}
  \begin{center}
  \includegraphics[width=0.31\textwidth]{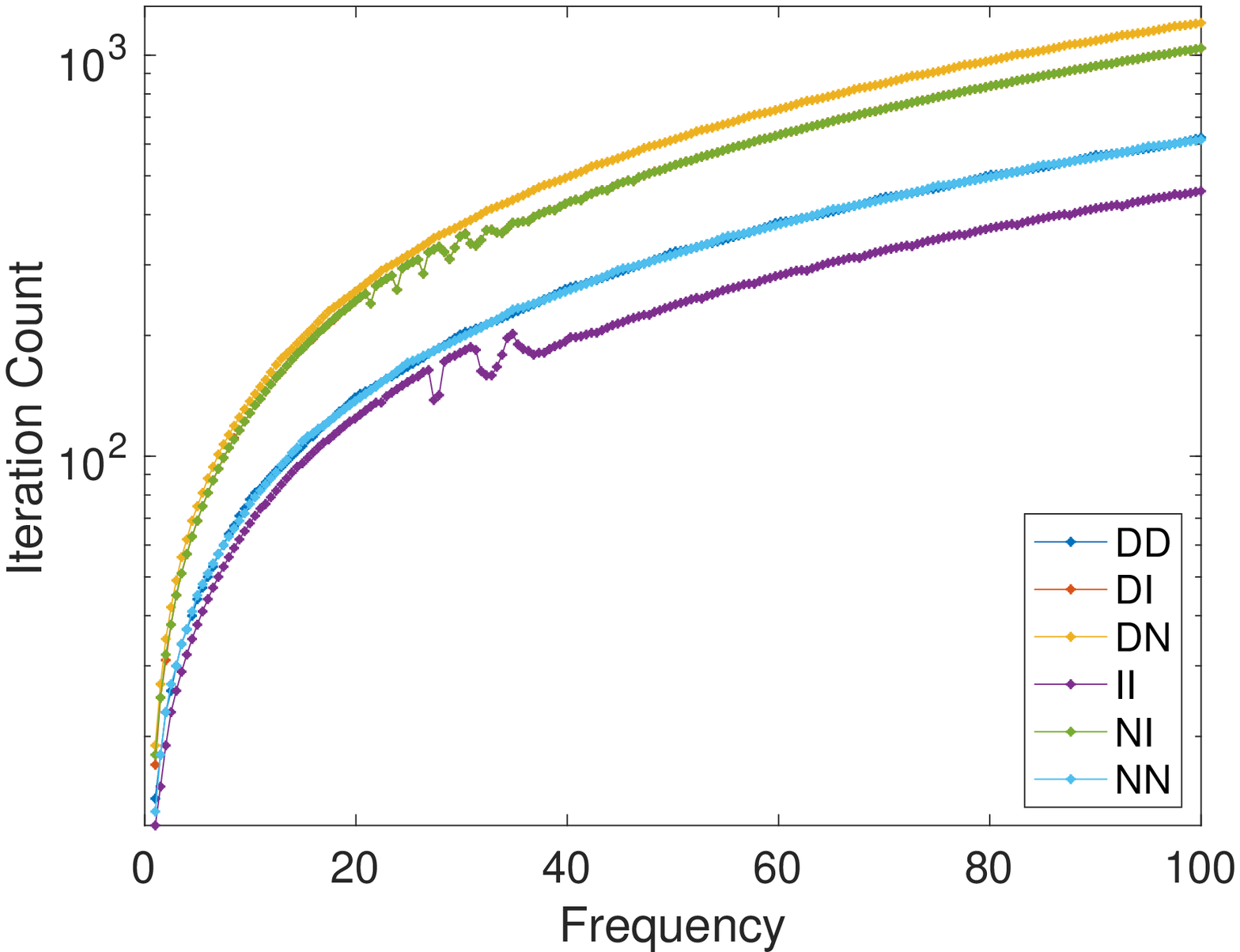}
  \includegraphics[width=0.32\textwidth]{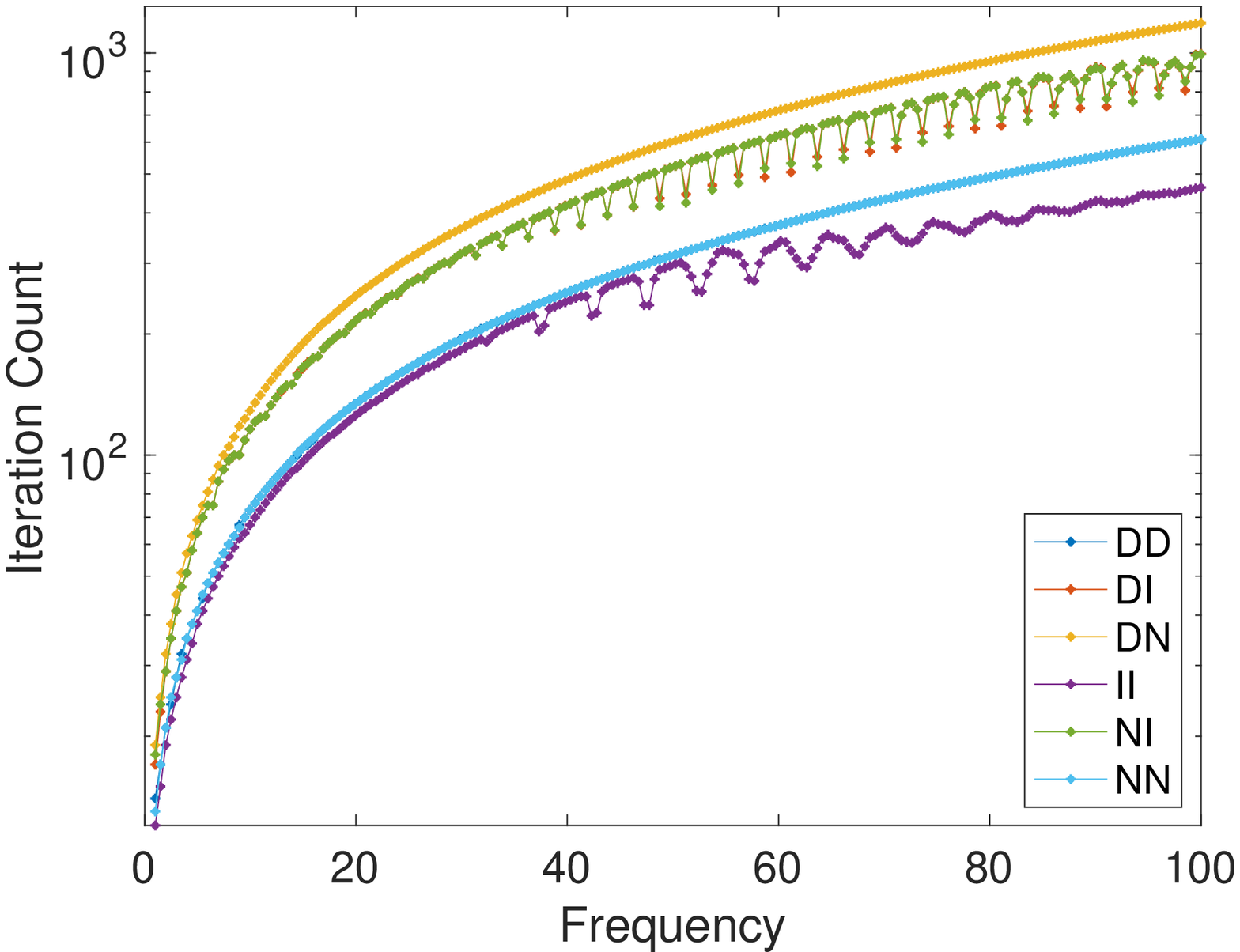}
  \includegraphics[width=0.32\textwidth]{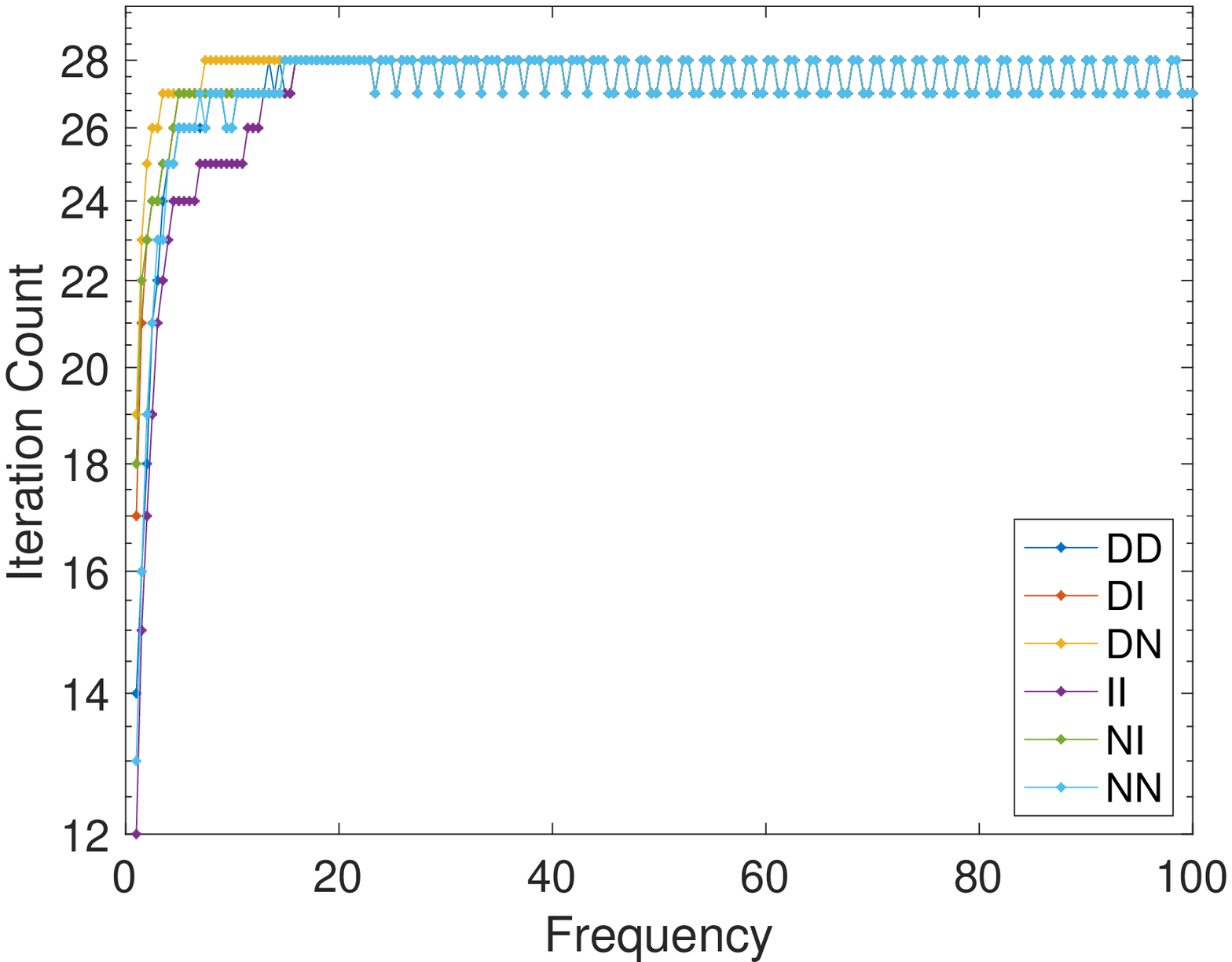}
  \caption{Number of iterations as a function of $\omega$ for different boundary conditions and damping parameters. Left to right: $\eta = 1/2\omega$, $1/2$, and $\omega/2$. Each legend entry indicates either Dirichlet (D), Neumann (N), or impedance (I) boundary conditions; the first letter
  at $x = -6$, and the second letter at $x = 6$. \label{fig:1d_iter_vs_freq}}
  \end{center}
\end{figure}

We report the number of iterations it takes to reach a GMRES residual smaller than $10^{-10}$ for the six possible combinations of Dirichlet, Neumann and impedance boundary conditions for 200 frequencies distributed evenly from 1 to 100. The results for three levels of damping are displayed in Figure \ref{fig:1d_iter_vs_freq}. On the left and middle of Figure \ref{fig:1d_iter_vs_freq} are damping parameters of $1/2\omega$ and $1/2$ respectively, from which it is clear that the scaling is sub-linear with increasing frequency. On the right in Figure \ref{fig:1d_iter_vs_freq} are results from a damping parameter that grows with frequency, $\omega/2$, which demonstrates a number of iterations that is both frequency independent and modest for a given GMRES tolerance. Interestingly, in this case the curve for each set of boundary conditions collapses to the same curve so that the iteration is insensitive to boundary conditions for a sufficiently large damping parameter.
\begin{remark}
As seen in the prequel \cite{WaveHoltz}, the impedance-impedance conditions take the fewest iterations to reach convergence for lower levels of damping. We point out the preceeding analysis assumes energy conserving boundary conditions to obtain estimates on the convergence rate of WaveHoltz as a fixed-point iteration. A different approach without the need for a Laplacian with a point-spectrum is needed to obtain rates depending on the specific boundary conditions.
\end{remark}
\subsection{Examples in Two Dimensions}
In this section we present experiments in two space dimensions.
For the following examples, we consider solving the Helmholtz equation
for the wedge model which we adapt from \cite{erlangga2006comparison,SoVHelm}. 
The domain is the rectangle $[0,600]\times[0,1000]$ with the (discontinuous)
speed of sound 
	\begin{align*}
		c(x) = 
		\begin{cases}
			c_1 = 2100, & y \le x/6 + 400, \\
			c_2 = 1000, & x/6 + 400 \le y \le 800 - x/3, \\ 
			c_3 = 2900, & \text{else}.
		\end{cases}
	\end{align*}
On the boundary of the rectangle we impose the
impedance boundary condition $w_t + c \nabla w \cdot \vec{n} = 0$.
For the spatial discretization we use the SIPDG method with a 
penalty parameter choice of $\gamma = (p+1)^2$, where $p=4$ 
is the polynomial order used in each element which results 
in a fifth order method. We discretize the domain with a total 
of 7680 triangular elements with a total of 115200 degrees of 
freedom. In time we use a fourth order Taylor method for time-stepping. For each example, 
we use the point-source 
	\begin{align*}
		f(x,y) = \omega^2 \delta(|x-x_0|)\delta(|y-y_0|),
	\end{align*}
where $x_0 = 300$, $y_0 = 0$, $\omega$ is the Helmholtz frequency,
and $\delta(z)$ is the usual Dirac delta function.
These examples were implemented in the MFEM finite element 
discretization library \cite{mfem}. 


\subsubsection{Convergence for Damped Helmholtz Equations}
We again study how the number of GMRES accelerated WHI iterations scale with the Helmholtz frequency $\omega$ for 
the exemplary wedge problem. 

We report the number of iterations it takes to reach a GMRES residual smaller than $10^{-10}$ 
for the frequencies $1, 2, \dots, 100$, with damping $\eta = \omega/2$ with either impedance or 
Neumann conditions on all sides of the rectangular domain.

\begin{figure}[htpb]
\graphicspath{{figures/}}
\begin{center}
\includegraphics[width=0.50\textwidth]{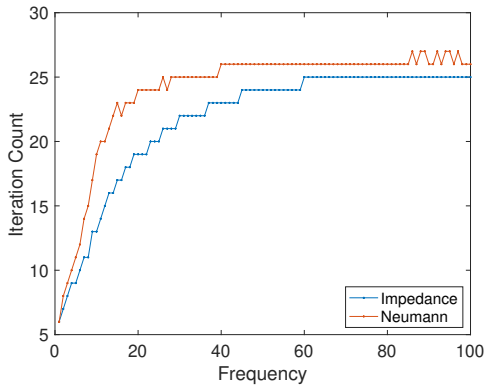}
\caption{Number of iterations to reach a GMRES tolerance of $10^{-10}$ for the wedge problem in 2D 
with all Neumann or all impedance boundary conditions.\label{fig:2d_ibc_conv}}
\end{center}
\end{figure}

The results for this experiment are shown in Figure \ref{fig:2d_ibc_conv}, 
from which it is clear that the number of iterations is essentially independent 
of frequency for larger frequencies as was the case in a single spatial dimension. 
We again note that energy conserving boundary conditions require more iterations than
the impedance case even in the presence of damping.

\begin{figure}[htpb]
\graphicspath{{figures/}}
\begin{center}
\includegraphics[width=0.49\textwidth]{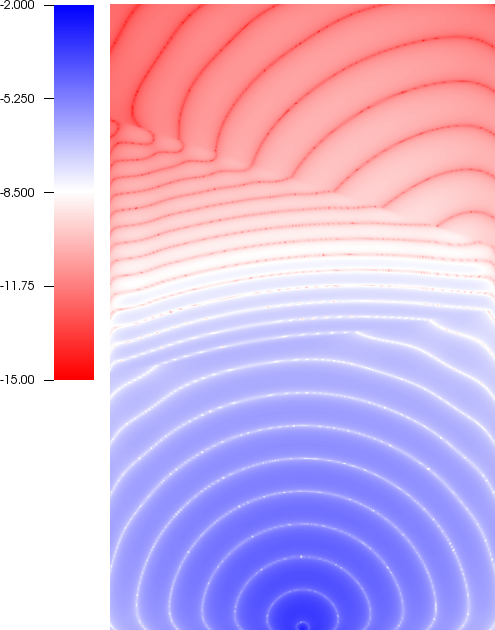}
\includegraphics[width=0.49\textwidth]{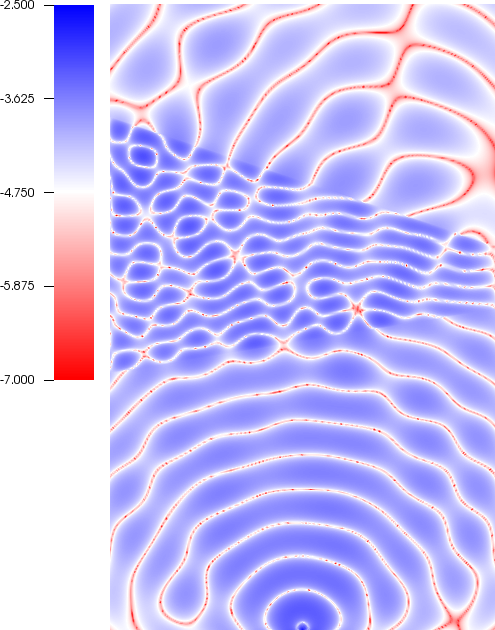}
\caption{In the above we plot the $\log_{10}$ of the absolute value of the real part 
of the Helmholtz solution with frequency $\omega = 40\pi$ for (Left) damping parameter $\eta = 20\pi$ and (Right) no damping. 
\label{fig:2d_freq_large}}
\end{center}
\end{figure}

For a final example, in Figure \ref{fig:2d_freq_large} we display the solution
of the damped (and undamped) Helmholtz equation using the GMRES accelerated WHI for 
a frequency of $\omega = 40\pi$ with damping $\eta = \omega/2$ and $0$, respectively.

\section{Summary and Future Work}
We have presented and extended analysis of the WaveHoltz iteration, 
an iterative method for solving the Helmholtz
equation, applied to wave equations with and without damping. 
The general iteration has the same rate of convergence as the energy conserving 
case presented in \cite{WaveHoltz}, but is a more general
and appropriate formulation for considering problems with impedance/Sommerfeld
boundary conditions. For problems with damping, the WaveHoltz iteration
always converges and numerical experiments verify the frequency independent
convergence of problems with sufficient levels of damping.

We have provided analysis of the interior impedance problem in a single dimension 
and constructed an example in which the worst-case convergence rate is realized, 
despite the numerical results of our previous paper indicating much more favorable 
scaling for non-energy conserving boundary conditions. We have additionally 
investigated higher order modified equation time-stepping schemes and shown that the 
WaveHoltz solution converges to the discrete Helmholtz solution to the order
matching the order of the chosen time-stepping scheme. In addition, we have presented
a method to \textit{completely} remove time-discretization error for centered
modified equation time-stepping schemes.

Finally, here we have only considered acoustic wave propagation. In future work 
we hope to apply the WaveHoltz iteration to elastic problems. Moreover, we have not
yet tried to leverage sweeping preconditioner ideas here
and hope to study the numerical and theoretical properties of these in the 
future.

\appendix 
\section{Proof of Lemma \ref{lemma::FilterBounds}}\label{sec:FilterBounds}
We show the results for the rescaled function
\begin{align*}
	\bar \gamma (r) := \gamma(r\omega) 
  &= \frac{2}{T} \int_0^T \left(\cos(\omega t) - \frac{1}{4}\right)\sin(r\omega t) \, dt
  \\
	&= \frac{1}{\pi} \int_0^{2\pi} \left(\cos(t) - \frac{1}{4}\right)\sin(r t)\, dt.
\end{align*}
By direct integration we get
	\begin{align*}
		\bar \gamma (r) 
    &= \frac{1}{\pi} \int_0^{2\pi} \frac{1}{2}\left(\sin((r+1) t) + \sin((r-1) t)\right) - \frac{1}{4}\sin(r t)\, dt 
    \\
    &= \frac{(1+3r^2) \sin^2(\pi r)}{2\pi r(r^2-1)} 
    \\
    &= \frac{\pi r(1+3r^2) \text{sinc}^2(r/2)}{2(r^2-1)},
	\end{align*}
where
	\begin{align*}
		\text{sinc}(r) = \frac{\sin(2\pi r)}{2\pi r}.
	\end{align*}
From \cite{WaveHoltz} we have the following expression for $\beta$:
	\begin{align*}
		\bar \beta (r) = \frac{1}{\pi} \int_0^{2\pi} \frac{1}{2}\left(\cos((r+1) t) + \cos((r-1) t)\right) - \frac{1}{4}\cos(r t)\, dt 
    &= \frac{(1+3r^2) \sin(2\pi r)}{4\pi r(r^2-1)} 
    \\
    &= \frac{(1+3r^2) \text{sinc}(r)}{2(r^2-1)}.
	\end{align*}
Then the eigenvalues of the WaveHoltz operator applied to the first order system
are
	\begin{align*}
		|\bar \mu(r)|^2 = \bar \beta^2 (r) + \bar \gamma^2 (r) = \frac{(1+3r^2)^2 \sin^2(\pi r)}{4\pi^2 r^2(r^2-1)^2} = \frac{(1+3r^2)^2 \text{sinc}^2(r/2)}{4(r^2-1)^2}.
	\end{align*}
We now first consider $0 \le r \le 0.5$ and note that $|\mu(r)|^2$ is a positive, increasing function
on this interval so that
	\begin{align*}
		|\bar \mu(r)|^2 \le |\bar \mu(1/2)| = \frac{49}{9\pi^2} \le 0.56.
	\end{align*}
For $1/2 \le r \le 3/2$ we instead center around $r = 1$ and get for $|\delta| \le 1/2$,
	\begin{align*}
		|\bar \mu(1+\delta)|^2 = \frac{(3(\delta+1)^2+1)^2 \sin^2(\pi \delta)}{4\pi^2 \delta^2(1+\delta)^2(2+\delta)^2} 
		=\frac{(3(\delta+1)^2+1)^2 \text{sinc}^2(\delta/2)}{4(1+\delta)^2(2+\delta)^2}.
	\end{align*}
We use the fact that $\sin(x) \le x - \tilde \alpha x^3$ in the interval $x \in [0, \pi]$ 
for any $\tilde \alpha \in [0,\pi^{-2}]$. This leads to the following estimate for the sinc
function
	\begin{align}\label{eqn::sincBound}
		0 \le \text{sinc}(r/2) \le 1-\alpha r^2, \quad r \in [-0.5, 0.5], \quad \alpha \in [0, 1].
	\end{align}
Using \eqref{eqn::sincBound} with $\alpha = 1$, gives 
	\begin{align*}
		|\bar \mu(1+\delta)|^2 \le \frac{(3(\delta+1)^2+1)^2 (1-\delta^2)^2}{4(1+\delta)^2(2+\delta)^2}
		&= \frac{(4 + 2\delta - 3\delta^2 - 3\delta^3)^2}{4(2+\delta)^2}
    \\
		&\le \frac{(4 + 2\delta - 2\delta^2 - \delta^3)^2}{4(2+\delta)^2} \\
		& = \left(1 - \frac{\delta^2}{2}\right)^2 \\
		& = 1 - \delta^2 + \frac{\delta^4}{4}
		\le 1 - \frac{15}{16}\delta ^2,
	\end{align*}
since $|\delta|<1/2$. A Taylor expansion around $\delta = 0$ for $|\delta| \le 1/2$
immediately gives the bound
	\begin{align*}
		\sqrt{1-\delta^2} \le 1 - \frac{\delta^2}{2} \implies |\mu(1+\delta)| 
		\le \sqrt{1-\frac{15}{16}\delta^2} \le 1 - \frac{15\delta^2}{32}.
	\end{align*}
If we consider $r \ge 3/2$, 
\begin{align*}
	|\bar \mu(r)|^2 = \frac{(1+3r^2)^2 \text{sinc}^2(r/2)}{4(r^2-1)^2} 
			   \le \frac{(1+3r^2)^2}{4(r^2-1)^2},
\end{align*}
which is a positive and decreasing function. It follows that
\begin{align*}
	|\bar \mu(r)|^2 \le \frac{(1+3(3/2)^2)^2}{4((3/2)^2-1)^2} \le 0.44,
\end{align*}
for $r\ge 3/2$. Finally, for a more general bound 
for $r > 1$ we have $1/(r+1) - 1/2r \ge 0$ so that 
\begin{align*}
	|\bar \mu(r)|^2 = \frac{(1+3r^2)^2 \sin^2(\pi r)}{4\pi^2 r^2(r^2-1)^2} 
			   \le 
			   \frac{(1+3r^2)^2}{4\pi^2 r^2(r^2-1)^2} 
			   &=
			   \frac{1}{\pi^2}
			   \left(\frac{1}{r+1} + \frac{1}{r-1} - \frac{1}{2r}\right)^2
         \\
			   &\le 
			   \left(\frac{3}{2\pi(r-1)}\right)^2,
\end{align*}
which gives
	\begin{align*}
		|\bar \mu(r)| \le \frac{3}{2\pi(r-1)}.
	\end{align*}
To prove \eqref{eqn::rhoest},
we use a Taylor expansion of $\bar \mu(r)$ about $r = 1$ in the interval $r \in (1/2,3/2)$,
	\begin{align*}
		|\bar \mu(1+\delta)| = 1 + \frac{\delta^2}{2}\frac{d^2}{dr^2}\left[|\bar \mu(r)|\right]_{r=1} + \frac{\delta^3}{6} 
	\bar R(\delta),
	\end{align*}
where $\bar R(\delta)$ is the remainder term. We note that 
by product rule we have
\begin{align*}
	\frac{d}{dr}|\bar\mu(r)| = \frac{1}{|\bar\mu|}(\bar\beta\bar\beta' + \bar\gamma\bar\gamma').
\end{align*}
Since
	\begin{align*}
		\frac{d}{dr}|\bar\mu(r)|^{-s} = -s |\bar \mu(r)|^{-s-1} \frac{d}{dr}|\bar\mu(r)|
		=
		\frac{-s}{|\bar \mu(r)|^{s+2}}(\bar\beta\bar\beta' + \bar\gamma\bar\gamma'),
	\end{align*}
by repeated product rule we can then show that 
	\begin{align*}
		\frac{d^3}{dr^3}|\bar\mu(r)| &= \frac{3}{|\bar\mu|^5}(\bar\beta\bar\beta' + \bar\gamma\bar\gamma')^2
		-\frac{1}{|\bar\mu|^3}(\bar\beta\bar\beta'' + (\bar\beta')^2 + \bar\gamma\bar\gamma'' + (\bar\gamma')^2)
		(1 +\bar\beta\bar\beta' + \bar\gamma\bar\gamma') \\
		& + \frac{1}{|\bar\mu|}(\bar\beta\bar\beta''' + 3\bar\beta'\bar\beta'' 
		+ \bar\gamma\bar\gamma''' + 3\bar\gamma'\bar\gamma'').
	\end{align*}
We note that $|\bar\mu(r)| \ge |\bar\mu(3/2)| \ge 1/\pi$ in the interval $1/2 \le r \le 3/2$,
and that we have the following bound
	\begin{align*}
		\sup_{r\geq 0} \left|\bar\beta^{(s)}(r)\right|\leq
	  \frac{1}{\pi}\int_0^{2\pi}t^s\left(1+\frac14\right)dt
	  =5\frac{2^{s-1}\,\pi^{s}}{s+1},
	\end{align*}
which similarly holds for $\sup_{r\geq 0} \left|\bar\gamma^{(s)}(r)\right|$
for $s = 0, 1, 2, \dots$. Thus by Taylor's theorem we have
	\begin{align*}
		|\bar R(\delta)|
		&\le 
		\sup_{1/2\le r \le 3/2} 
		\left|\frac{d^3}{dr^3}|\bar\mu(r)|\right| 
    \\
		& \le 
		\frac{3}{|\bar\mu(3/2)|^5}\frac{25^2 \pi^2}{4}
		+
		\frac{3}{|\bar\mu(3/2)|^3}\left(\frac{50 \pi^2}{3} + \frac{25 \pi^2}{2}\right)
		\left(1 + \frac{25 \pi}{2}\right)
		+
		\frac{3 \cdot 75 \pi^3}{|\bar\mu(3/2)|}\\
		& \le 
		\frac{3}{4} 25^2 \pi^7
		+
		3\pi^3 \left(\frac{50 \pi^2}{3} + \frac{25 \pi^2}{2}\right)
		\left(1 + \frac{25 \pi}{2}\right)
		+
		3 \cdot 75 \pi^4\\
		& = 
		\frac{25 \pi^4}{4}\left(36 + 20 \pi + 250 \pi^2 + 75 \pi^3\right).
	\end{align*}

Then, $|R(\delta)| \le |\bar R(\delta)|/6$.
Finally,
	\begin{align*}
		\frac{d^2}{dr^2}\left[|\bar \mu(r)|\right]_{r=1} = \frac16 (3 - 2 \pi^2) = - 2 b_1.
	\end{align*}
\section{Wave Equation Extension}\label{sec:WaveExtension}
Let $\Omega = (-\infty,0)$ and let $f\in L^2(\Omega)$ be compactly supported 
in $\Omega$ away from $x=0$. Additionally, assume 
$1/c^2 \in L^1_{\text{loc}}(\Omega)$ with $c(0) = c_0$ on the 
interval $[-\delta,0]$ for some $\delta >0$.
We consider the semi-infinite problem
  \begin{eqnarray*}
    &&  w_{tt} =   \frac{\partial}{\partial x} \left[c^2(x) \frac{\partial}{\partial x}w\right]
      - \text{Re}\{f(x) e^{i \omega t}\}, \quad x \le 0, \ \ 0 \le t \le T, \nonumber \\
    &&  w(0,x) = v_0(x), \quad w_t(0,x) = v_1(x), \nonumber \\
    &&  \alpha w_t(t,0) +\beta c_0 w_x(t,0)=0.
  \end{eqnarray*}
Let $\tilde w$ solve the extended wave equation
  \begin{eqnarray*}
    &&  \tilde w_{tt} =   \frac{\partial}{\partial x} \left[\tilde c^2(x) \frac{\partial}{\partial x}\tilde w\right]
      - \text{Re}\{\tilde f(x) e^{i \omega t}\}, \quad x \in \mathbb{R}, \ \ 0 \le t \le T, \nonumber \\
    &&  \tilde w(0,x) = \tilde v_0(x), \quad \tilde w_t(0,x) = \tilde v_1(x), \nonumber \\
    &&  \alpha \tilde w_t(t,0) +\beta c_0 \tilde w_x(t,0)=0,
  \end{eqnarray*}
where $\tilde f$ is a zero extension, $\tilde c$ is 
the extended wavespeed
  \begin{align*}
    \tilde c(x) = 
    \begin{cases}
      c_0, & -\delta < x\le0, \\
      \tilde c_0, & x >0,
    \end{cases}
  \end{align*}
and (\textbf{I}) $\tilde v_0 \in H^1(\mathbb{R})$ and $\tilde v_1 \in L^2(\mathbb{R})$
are extensions of the initial data.
We choose the extensions of $v_0$ and $v_1$ such 
that
  \begin{align*}
    \tilde v_1(x) + \tilde v_0'(x) = 0, \quad x>0, \quad (\textbf{II})
  \end{align*}
so that the wave solution in the region $x>0$
satisfies the condition $w_t + w_x = 0$ at $x=0$,
ensuring no data propagates into the original
domain $x<0$. In particular, we may take 
$\tilde v_0$ to be constant and $\tilde v_1 \equiv 0$.
Moreover, since $c$ is constant
in $[-\delta,0]$ the solution will then be of the form
  \begin{align*}
    \tilde w(t,x) = 
    \begin{cases}
      w_L(x+c_0t) + w_R(x-c_0 t), & -\delta \le x\le0, \\
      w_T(x-\tilde c_0 t), & x >0,
    \end{cases}
  \end{align*}
for some functions $w_L, w_R$, and $w_T$.
At $x=0$ where $\tilde c$ is (potentially) discontinuous,
the weak solution satisfies the interface conditions that 
$\tilde w$ and $\tilde c^2 \tilde w_x$ are both 
continuous. These requirements lead to the 
relations
  \begin{align*}
    w_L(c_0t) + w_R(-c_0t) = w_T(-\tilde c_0 t), \\ 
    c_0^2(w_L'(c_0t) + w_R'(-c_0 t)) = \tilde c_0^2 w_T'(-c_0t).
  \end{align*}
It follows that 
  \begin{align*}
    \tilde w_t(t,0^-) = c_0 (w_L'(c_0t) - w_R'(-c_0t)) = -\tilde c_0w_T'(-\tilde c_0 t), 
    \quad c_0\tilde w_x(t,0^-) = \frac{\tilde c_0^2}{c_0} w_T'(-c_0t),
  \end{align*}
so that the impedance condition
  \begin{align*}
    \alpha \tilde w_t(t,0^-) + \beta c_0 \tilde w_x(t,0^-) = \left(-\alpha \tilde c_0 + \beta \frac{\tilde c_0^2}{c_0}\right)w_T'(-c_0t) = 0,
  \end{align*}
is satisfied if 
  \begin{align*}
    \tilde c_0 = \frac{\alpha}{\beta} c_0. \quad (\textbf{III})
  \end{align*}
With this choice of the extended wavespeed $\tilde c_0$, both 
$\tilde w$ and $w$ satisfy the same PDE and condition at $x=0$ so 
that they must be equal for $x<0$. In summary,
if we have that conditions (\textbf{I}-\textbf{III}) are satisfied,
we have that $\tilde w(t,x) = w(t,x)$ for $x<0$. We note that
a similar argument can be made for an interior impedance problem
on a bounded domain, $a \le x \le b$, to 
a problem on $\mathbb{R}$. In this case, assuming
$c(a) = c_a$, $c(b) = c_b$ where $c$ is constant near the endpoints,
then the following problem has $\tilde w(t,x) = w(t,x)$ for 
$a \le x \le b$:
  \begin{eqnarray*}
    &&  \tilde w_{tt} =   \frac{\partial}{\partial x} \left[\tilde{c}^2(x) \frac{\partial}{\partial x}\tilde w\right]
      - \text{Re}\{\tilde f(x) e^{-i\omega t}\}, \quad x\in \mathbb{R}, \ \ 0 \le t \le T, \nonumber \\
    &&  \tilde w(0,x) = \tilde{v}_0(x), \quad \tilde w_t(0,x) = \tilde{v}_1(x), \nonumber
  \end{eqnarray*}
where $\tilde{v}_0$ and $\tilde c$ are the constant extensions
(with $\gamma = \alpha/\beta$)
  \begin{align*}
    \tilde{v}_0(x) = 
      \begin{cases}
        v_0(a_0), & x < a, \\
        v_0(x), & a \le x \le b, \\
        v_0(b_0), & b < x ,
      \end{cases}
      \quad 
      \tilde{c}(x) = 
        \begin{cases}
          \gamma c_a, & x < a, \\
          c(x), & a \le x \le b, \\
          \gamma c_b, & b < x,
        \end{cases}
    \end{align*}
and $\tilde{v}_1$, $\tilde f$ are zero extensions of $v_1$ and $f$, respectively.

Since the solutions to the wave equation have finite speed of 
propagation, we may replace the domain $\mathbb{R}$ for $\tilde w$ by 
a large enough domain with any boundary condition given that 
any reflections at the new boundary do not re-enter the region $a \le x \le b$. 
Let $\tilde a < a - c_aT/2$ and $\tilde b > b + c_bT/2$.
We define the extension operator $E$ such that $[v_0,v_1]^T \rightarrow [\tilde v_0,\tilde v_1]^T$
where $\tilde v_0$ and $\tilde c$ are the extensions as above
and $\tilde{v}_1$, $\tilde f$ are zero extensions of $v_1$ and $f$, respectively.
We now consider the (finite interval) extended problem with homogeneous Neumann 
conditions
  \begin{eqnarray*}
    &&  \tilde w_{tt} =   \frac{\partial}{\partial x} \left[\tilde{c}^2(x) \frac{\partial}{\partial x}\tilde w\right]
      - \text{Re}\{\tilde f(x) e^{-i\omega t}\}, \quad \tilde a\le x \le \tilde b, \ \ 0 \le t \le T, \nonumber \\
    &&  \tilde w(0,x) = \tilde{v}_0(x), \quad \tilde w_t(0,x) = \tilde{v}_1(x), \nonumber \\
    &&  \tilde w_x(t,\tilde a) = 0, \quad \tilde w_x(t,\tilde b)=0.
  \end{eqnarray*}
Defining the projection operator $P$ as the restriction 
of $\tilde w$ to $a\le x\le b$ then it follows that $P \tilde w = w$ 
where $w$ is the original wave solution to the interior impedance problem.

\section{Well-definedness of modified frequencies}\label{sec:well-defined}

Here we show that the modified frequency $\tilde\omega$
is well-defined. 
This is given by the following lemma.
\begin{lemma}\label{lem::welldefined}
For each $\omega\in\Real$ satisfying $0<\Delta t\omega\leq 1$
there is a modified frequency $\tilde{\omega}\in\Real$
which is the smallest positive real number satisfying
\begin{align}\label{eqn::freq_relation}
      \sin^2(\omega\Delta t/2) = \sum_{j=1}^m \frac{(-1)^{j+1}\left(\Delta t \tilde \omega\right)^{2j}}{2(2j)!}.
\end{align}
Moreover, there is a constant $C_m<1$ that only depends on $m$
such that
\begin{align}\label{eqn::omegaerr}
0<\Delta t\tilde\omega\leq 2, \qquad
|\omega-\tilde{\omega}|\leq C_m \Delta t^{2m}\omega^{2m+1},\qquad 
C_m := 5/(2m+2)!
\end{align}
and for all $0\leq \Delta t \lambda \leq 2$, it holds that
\begin{align}\label{eqn::lambound}
\left|\sum_{j=1}^m \frac{(-1)^{j+1}\left(\Delta t \lambda\right)^{2j}}{2(2j)!}\right|\leq 1.
\end{align}

\end{lemma}
\begin{proof}
We define the polynomial
  \begin{align}\label{eqn::p_func}
    p(x) = \sum_{j=1}^m \frac{(-1)^{j+1}x^j}{2(2j)!} - \sin^2(\omega\Delta t/2),
  \end{align}
and note that $(\Delta t \tilde \omega)^2$ is a root of $p(x)$.
On the interval $[0,1/2]$ we have that $\sin^2(x)$ is increasing so that 
  \begin{align*}
    0 < \sin^2(\omega\Delta t/2) \le \sin^2(1/2) \le 0.23 < 1,
  \end{align*}
immediately giving $p(0) < 0$. Moreover,
  \begin{align}\label{eq::p4def}
    p(4) 
    &= -\frac{1}{2}\sum_{j=1}^m \frac{(-1)^{j}2^{2j}}{(2j)!} - \sin^2(\omega\Delta t/2) \nonumber
    \\
    & =  -\frac{1}{2}\sum_{j=1}^\infty \frac{(-1)^{j}2^{2j}}{(2j)!} + \frac{1}{2}\sum_{j=m+1}^\infty \frac{(-1)^{j}2^{2j}}{(2j)!} - \sin^2(\omega\Delta t/2) \nonumber\\
    & =  \sin^2(1) + \frac{1}{2}\sum_{j=m+1}^\infty \frac{(-1)^{j}2^{2j}}{(2j)!} - \sin^2(\omega\Delta t/2).
  \end{align}
%
%
%
%
We note that, 
  \begin{align*}
    \sum_{j=2}^\infty \frac{2^{2j}}{(2j)!} = \sum_{j=0}^\infty \frac{2^{2j}}{(2j)!} - 3 = \cosh(2) - 3,
  \end{align*}
so that 
  \begin{align}\label{eq::coshest}
    \left|\frac{1}{2}\sum_{j=m+1}^\infty \frac{(-1)^{j}2^{2j}}{(2j)!}\right| 
    < \frac{1}{2}\sum_{j=m+1}^\infty \frac{2^{2j}}{(2j)!} 
    \leq\frac{1}{2}\sum_{j=2}^\infty \frac{2^{2j}}{(2j)!}
    =
    \frac{1}{2}(\cosh(2) - 3).
  \end{align}
  Since also $\sin^2(\omega\Delta t/2)\leq \sin^2(1/2)$,
we get
  \begin{align*}
    p(4)> 
    \sin^2(1) -\frac{1}{2}(\cosh(2) - 3) - \sin^2(1/2)
    \approx 0.097
     > 0.
  \end{align*}
By the intermediate value theorem, it follows that $p(x)$ has a root in the interval $(0,4)$.
We next need to show that $p'(x) \ne 0$ on this interval to guarantee the root is unique.
Taking a derivative,
  \begin{align*}
    \frac{d}{dx} p(x) = \sum_{j=1}^m \frac{(-1)^{j+1}j x^{j-1}}{(2j)!}
    =
    \frac{1}{2}\left[1 + \sum_{j=2}^m \frac{(-1)^{j+1}x^{j-1}}{(2j-1)!}\right].
  \end{align*}
We then have
  \begin{align*}
    \sum_{j=2}^m \frac{(-1)^{j+1}x^{j-1}}{(2j-1)!}
    \ge 
    -\sum_{j=2}^m \frac{2^{j-1}}{(2j-1)!}
    &=
    -\frac{1}{\sqrt{2}}
    \sum_{j=2}^m \frac{\sqrt{2}^{2j-1}}{(2j-1)!}
    \\
    &\ge
    -\frac{1}{\sqrt{2}} \left[\sum_{j=1}^\infty \frac{\sqrt{2}^{2j-1}}{(2j-1)!} -\sqrt{2}\right]
    \\
    &=
    -\frac{\sinh(\sqrt{2})}{\sqrt{2}}+1,
  \end{align*}
so that 
  \begin{align}\label{eqn::p_bound}
    \frac{d}{dx} p(x) = \frac{1}{2}\left[1 + \sum_{j=2}^m \frac{(-1)^{j+1}x^{j-1}}{(2j-1)!}\right]
    \ge 
    \frac{1}{2} -\frac{\sinh(\sqrt{2})}{\sqrt{2}} + 1> \frac{1}{10} > 0.
  \end{align}
This gives that there is a unique, positive, real-valued $\tilde \omega$ with $0<\Delta t \tilde \omega < 2$
that satisfies the relation \eqref{eqn::freq_relation}, 
showing the first part of the lemma.
For the last part
we let $x = \Delta t \omega\in(0,1]$, $\tilde x = \Delta t \tilde \omega\in(0,2)$,
and $R_m = p(x^2) - p(\tilde x^2)$ where $p(x)$ is defined as in \eqref{eqn::p_func}. 
By the mean value theorem we have
  \begin{align*}
    |R_m| = |p(x^2) - p(\tilde x^2)| = |(x^2 - \tilde x^2) p'(\xi)| =  |x - \tilde x| |x + \tilde x| |p'(\xi)|,
  \end{align*}
for some $\xi \in (0,2)$, so that
  \begin{align*}
    |x - \tilde x| \le \frac{|R_m|}{|x + \tilde x| |p'(\xi)|}.
  \end{align*}
Since $\tilde x^2$ is a root of $p$, a Taylor series estimate gives
  \begin{align*}
    |R_m| = |p(x^2)|
          = \left| \sum_{j=1}^m \frac{(-1)^{j+1}x^j}{2(2j)!} - \sin^2(\omega\Delta t/2) \right|
          &= \left|  \sum_{j=m+1}^\infty \frac{(-1)^{j+2}\left(\Delta t \omega\right)^{2j}}{2(2j)!} \right|
          \\ &\le \frac{\left(\Delta t \omega\right)^{2m+2}}{2(2m+2)!},
  \end{align*}
which gives
  \begin{align*}
    |x - \tilde x| \le \frac{|R_m|}{|x + \tilde x| |p'(\xi)|} \le \frac{\Delta t^{2m+1} \omega^{2m+2}}{2(2m+2)!(\omega + \tilde \omega) |p'(\xi)|}
    \\\implies 
    |\omega - \tilde \omega| \le \frac{\Delta t^{2m} \omega^{2m+2}}{2(2m+2)!(\omega + \tilde \omega) |p'(\xi)|}.
  \end{align*}
By \eqref{eqn::p_bound} we have that $|p'(x)| > 1/10$ in $[0,4]$, which finally gives
  \begin{align}\label{eqn::Freq_error}
    |\omega - \tilde \omega| \le \frac{10 \Delta t^{2m} \omega^{2m+2}}{2(2m+2)!(\omega + \tilde \omega)} 
    \le  \frac{5 \Delta t^{2m} \omega^{2m+1}}{(2m+2)!
    }
    =: C_m\Delta t^{2m} \omega^{2m+1}.
  \end{align}
  We finally prove 
  \eqref{eqn::lambound}.
  It is trivially true for $\lambda=0$ so we assume
  that $\lambda>0$.
%
%
We define $a_k = (\Delta t \lambda)^{2k}/2(2k)!$, and note that
$a_{k+1}/a_{k} < 1$ so that $a_k > a_{k+1}>0$ and 
that $a_1 =  (\Delta t \lambda)^{2}/4 \le 1$. 
Letting
  \begin{align*}
    \tilde m = 
    \begin{cases}
      m, & m\text{ odd},\\
      m+1, & m\text{ even},
    \end{cases}
  \end{align*}
then we have
  \begin{align*}\label{eqn::lam_welldefined}
    \sum_{k=1}^m \frac{(-1)^{k+1}\left(\Delta t \lambda\right)^{2k}}{2(2k)!}
    =
    \sum_{k=1}^m (-1)^{k+1}a_k
    &\le 
    \sum_{k=1}^{\tilde m} (-1)^{k+1}a_k
    \\
    &=
    a_1
    -
    \sum_{k=1}^{(\tilde m-1)/2}
    (a_{2k}-a_{2k+1})
    \le 
    a_1 
    \le 1,
  \end{align*}
If instead
  \begin{align*}
    \tilde m = 
    \begin{cases}
      m+1, & m\text{ odd},\\
      m, & m\text{ even},
    \end{cases}
  \end{align*}
then we have the bound
  \begin{align*}
    \sum_{k=1}^m \frac{(-1)^{k+1}\left(\Delta t \lambda\right)^{2k}}{2(2k)!}
    =
    \sum_{k=1}^m (-1)^{k+1}a_k
    \ge 
    \sum_{k=1}^{\tilde m} (-1)^{k+1}a_k
    =
    \sum_{k=1}^{\tilde m/2}
    (a_{2k-1}-a_{2k})
    >
    0,
  \end{align*}
proving 
\eqref{eqn::lambound}. 
\end{proof}

\bibliographystyle{amsplain}
\bibliography{garcia}
\end{document}